
\catcode'32=9
\magnification=1200


\voffset=1cm
\hoffset=.5cm
\font\tenpc=cmcsc10

\font\eightrm=cmr8
\font\eighti=cmmi8
\font\eightsy=cmsy8
\font\eightbf=cmbx8
\font\eighttt=cmtt8
\font\eightit=cmti8
\font\eightsl=cmsl8
\font\sixrm=cmr6
\font\sixi=cmmi6
\font\sixsy=cmsy6
\font\sixbf=cmbx6

\skewchar\eighti='177 \skewchar\sixi='177
\skewchar\eightsy='60 \skewchar\sixsy='60


\font\tengoth=eufm10
\font\tenbboard=msbm10
\font\eightgoth=eufm10 at 8pt
\font\eightbboard=msbm8
\font\sevengoth=eufm7
\font\sevenbboard=msbm7
\font\sixgoth=eufm6
\font\fivegoth=eufm5

\newfam\gothfam
\newfam\bboardfam

\catcode`\@=11

\def\raggedbottom{\topskip 10pt plus 36pt
\r@ggedbottomtrue}
\def\pc#1#2|{{\bigf@ntpc #1\penalty
\@MM\hskip\z@skip\smallf@ntpc #2}}

\def\tenpoint{%
  \textfont0=\tenrm \scriptfont0=\sevenrm \scriptscriptfont0=\fiverm
  \def\rm{\fam\z@\tenrm}%
  \textfont1=\teni \scriptfont1=\seveni \scriptscriptfont1=\fivei
  \def\oldstyle{\fam\@ne\teni}%
  \textfont2=\tensy \scriptfont2=\sevensy \scriptscriptfont2=\fivesy
  \textfont\gothfam=\tengoth \scriptfont\gothfam=\sevengoth
  \scriptscriptfont\gothfam=\fivegoth
  \def\goth{\fam\gothfam\tengoth}%
  \textfont\bboardfam=\tenbboard \scriptfont\bboardfam=\sevenbboard
  \scriptscriptfont\bboardfam=\sevenbboard
  \def\bboard{\fam\bboardfam}%
  \textfont\itfam=\tenit
  \def\it{\fam\itfam\tenit}%
  \textfont\slfam=\tensl
  \def\sl{\fam\slfam\tensl}%
  \textfont\bffam=\tenbf \scriptfont\bffam=\sevenbf
  \scriptscriptfont\bffam=\fivebf
  \def\bf{\fam\bffam\tenbf}%
  \textfont\ttfam=\tentt
  \def\tt{\fam\ttfam\tentt}%
  \abovedisplayskip=12pt plus 3pt minus 9pt
  \abovedisplayshortskip=0pt plus 3pt
  \belowdisplayskip=12pt plus 3pt minus 9pt
  \belowdisplayshortskip=7pt plus 3pt minus 4pt
  \smallskipamount=3pt plus 1pt minus 1pt
  \medskipamount=6pt plus 2pt minus 2pt
  \bigskipamount=12pt plus 4pt minus 4pt
  \normalbaselineskip=12pt
  \setbox\strutbox=\hbox{\vrule height8.5pt depth3.5pt width0pt}%
  \let\bigf@ntpc=\tenrm \let\smallf@ntpc=\sevenrm
  \let\petcap=\tenpc
  \normalbaselines\rm}
\def\eightpoint{%
  \textfont0=\eightrm \scriptfont0=\sixrm \scriptscriptfont0=\fiverm
  \def\rm{\fam\z@\eightrm}%
  \textfont1=\eighti \scriptfont1=\sixi \scriptscriptfont1=\fivei
  \def\oldstyle{\fam\@ne\eighti}%
  \textfont2=\eightsy \scriptfont2=\sixsy \scriptscriptfont2=\fivesy
  \textfont\gothfam=\eightgoth \scriptfont\gothfam=\sixgoth
  \scriptscriptfont\gothfam=\fivegoth
  \def\goth{\fam\gothfam\eightgoth}%
  \textfont\bboardfam=\eightbboard \scriptfont\bboardfam=\sevenbboard
  \scriptscriptfont\bboardfam=\sevenbboard
  \def\bboard{\fam\bboardfam}%
  \textfont\itfam=\eightit
  \def\it{\fam\itfam\eightit}%
  \textfont\slfam=\eightsl
  \def\sl{\fam\slfam\eightsl}%
  \textfont\bffam=\eightbf \scriptfont\bffam=\sixbf
  \scriptscriptfont\bffam=\fivebf
  \def\bf{\fam\bffam\eightbf}%
  \textfont\ttfam=\eighttt
  \def\tt{\fam\ttfam\eighttt}%
  \abovedisplayskip=9pt plus 2pt minus 6pt
  \abovedisplayshortskip=0pt plus 2pt
  \belowdisplayskip=9pt plus 2pt minus 6pt
  \belowdisplayshortskip=5pt plus 2pt minus 3pt
  \smallskipamount=2pt plus 1pt minus 1pt
  \medskipamount=4pt plus 2pt minus 1pt
  \bigskipamount=9pt plus 3pt minus 3pt
  \normalbaselineskip=9pt
  \setbox\strutbox=\hbox{\vrule height7pt depth2pt width0pt}%
  \let\bigf@ntpc=\eightrm \let\smallf@ntpc=\sixrm
  \normalbaselines\rm}

\tenpoint

\frenchspacing


\newif\ifpagetitre
\newtoks\auteurcourant \auteurcourant={\hfil}
\newtoks\titrecourant \titrecourant={\hfil}

\def\appeln@te{}
\def\vfootnote#1{\def\@parameter{#1}\insert\footins\bgroup\eightpoint
  \interlinepenalty\interfootnotelinepenalty
  \splittopskip\ht\strutbox 
  \splitmaxdepth\dp\strutbox \floatingpenalty\@MM
  \leftskip\z@skip \rightskip\z@skip
  \ifx\appeln@te\@parameter\indent \else{\noindent #1\ }\fi
  \footstrut\futurelet\next\fo@t}

\pretolerance=500 \tolerance=1000 \brokenpenalty=5000
\newdimen\hmargehaute \hmargehaute=0cm
\newdimen\lpage \lpage=13.3cm
\newdimen\hpage \hpage=20cm
\newdimen\lmargeext \lmargeext=1cm
\hsize=12cm
\vsize=18cm
\parskip 0pt
\parindent=12pt

\def\margehaute{\vbox to \hmargehaute{\vss}}%
\def\margebasse{\vss}

\output{\shipout\vbox to \hpage{\margehaute\nointerlineskip
  \corpsdepage\margebasse}
  \advancepageno \global\pagetitrefalse
  \ifnum\outputpenalty>-20000 \else\dosupereject\fi}

\def\corpsdepage{\hbox to \lpage{\hss\pagetexte\hskip\lmargeext}}
\def\pagetexte{\vbox{\makeheadline\pagebody\makefootline}}
\headline={\ifpagetitre\titleheadline \else
  \ifodd\pageno\rightheadline \else\leftheadline\fi\fi}
\def\leftheadline{\eightpoint\hfil\the\auteurcourant\hfil}
\def\rightheadline{\eightpoint\hfil\the\titrecourant\hfil}
\def\titleheadline{\hfill}
\pagetitretrue

\def\footnoterule{\kern-6\p@
  \hrule width 2truein \kern 5.6\p@} 

\def\pd#1#2 {\pc#1#2| }

\def\pointir{\discretionary{.}{}{.\kern.35em---\kern.7em}\nobreak
\hskip 0em plus .3em minus .4em }

\def\abstract#1{\vbox{\eightpoint \pc ABSTRACT|\pointir #1}}

\def\titre#1|{\message{#1}
              \par\vskip 30pt plus 24pt minus 3pt\penalty -1000
              \vskip 0pt plus -24pt minus 3pt\penalty -1000
              \centerline{\bf #1}
              \vskip 5pt
              \penalty 10000 }

\def\section#1|{\par\vskip .3cm
                {\bf #1}\pointir}

\def\ssection#1|{\par\vskip .2cm
                {\it #1}\pointir}

\long\def\th#1|#2\finth{\par\medskip
              {\petcap #1\pointir}{\it #2}\par\smallskip}

\long\def\tha#1|#2\fintha{\par\medskip
                    {\petcap #1.}\par\nobreak{\it #2}\par\smallskip}

\def\rem#1|{\par\medskip
            {{\it #1}\pointir}}

\def\rema#1|{\par\medskip
             {{\it #1.}\par\nobreak }}

\def\ieme{\raise 1ex\hbox{\pc{}i\`eme|}}
\def\omini{\raise 1ex\hbox{\pc{}o|}}
\def\emini{\raise 1ex\hbox{\pc{}e|}}
\def\ermini{\raise 1ex\hbox{\pc{}er|}}
\def\remini{\raise 1ex\hbox{\pc{}re|}}

\def\article#1|#2|#3|#4|#5|#6|#7|
    {{\leftskip=7mm\noindent
     \hangindent=2mm\hangafter=1
     \llap{[#1]\hskip.35em}{#2}\pointir
     #3, {\sl #4}, t.\nobreak\ {\bf #5}, {\oldstyle #6},
     p.\nobreak\ #7.\par}}
\def\livre#1|#2|#3|#4|
    {{\leftskip=7mm\noindent
    \hangindent=2mm\hangafter=1
    \llap{[#1]\hskip.35em}{#2}\pointir
    {\sl #3}\pointir #4.\par}}
\def\divers#1|#2|#3|
    {{\leftskip=7mm\noindent
    \hangindent=2mm\hangafter=1
     \llap{[#1]\hskip.35em}{#2}\pointir
     #3.\par}}
\mathchardef\conj="0365

\def\qed{\quad\raise -2pt\hbox{\vrule\vbox to 10pt{\hrule width 4pt
\vfill\hrule}\vrule}}

\def\virg{\raise 2pt\hbox{,}}   

\long\def\entourer#1{\hbox{\vrule\vbox{\hrule\hbox{\kern15pt\vbox{\kern5pt
{#1}\kern5pt}\kern15pt}\hrule}\vrule}}
\def\\S {\vskip 5pt\hskip .5cm plus .1cm minus .1cm\relax}

\def\enonce#1|#2\finenonce{{\par\leftskip=36pt
\noindent\hbox to 0pt{\kern-\leftskip#1\hfill}{#2}\par}}

\def\senonce#1|#2\finsenonce{{\par\leftskip=36pt
\noindent\hbox to 0pt{\kern-24pt #1\hfill}{#2}\par}}

\def\ssenonce#1|#2\finssenonce{{\par\leftskip=48pt
\noindent\hbox to 0pt{\kern-24pt #1\hfill}{#2}\par}}

\def\decale#1|{\par\noindent\hskip 28pt\llap{#1}\kern 5pt}
\def\decaledecale#1|{\par\noindent\hskip 34pt\llap{#1}\kern 5pt}
\def\titrea#1|#2|{\message{#1 #2}
  \par\vskip.5cm plus .1cm minus .1cm\penalty -1000
  \centerline{\bf #1}
  \centerline{\bf #2}
  \vskip 5pt
  \penalty 10000 }
\def\sectiona#1|{\par\vskip .3cm
  {\bf #1.}
  \par\nobreak\vskip 3pt }
\def\ssectiona#1|{\par\vskip .2cm
  {\it #1.}
  \par\nobreak\vskip 2pt }

\def\rest#1{\ifinner\setbox1=\hbox{$\textstyle{#1}$}
            \else\setbox1=\hbox{$\displaystyle{#1}$}\fi
            \dimen1=\ht1
            \advance\dimen1 by\dp1
            \divide\dimen1 by 2
            \box1\lower 2pt\hbox{$\left|\vbox to\dimen1{}\right.$}}

\def\displaylinesno#1{\displ@y\halign{
\hbox to\displaywidth{$\@lign\hfil\displaystyle##\hfil$}&
\llap{$##$}\crcr#1\crcr}}

\def\ldisplaylinesno#1{\displ@y\halign{ 
\hbox to\displaywidth{$\@lign\hfil\displaystyle##\hfil$}&
\kern-\displaywidth\rlap{$##$}\tabskip\displaywidth\crcr#1\crcr}}

\def\Eqalign#1{\null\,\vcenter{\openup\jot\m@th\ialign{
\strut\hfil$\displaystyle{##}$&$\displaystyle{{}##}$\hfil
&&\quad\strut\hfil$\displaystyle{##}$&$\displaystyle{{}##}$\hfil
\crcr#1\crcr}}\,}


\def\matrixd#1{\null \,\vcenter {\normalbaselines \m@th
\ialign {\hfil $##$&&\quad \hfil $##$\crcr
\mathstrut \crcr \noalign {\kern -\baselineskip } #1\crcr
\mathstrut \crcr \noalign {\kern -\baselineskip }}}\,}

\def\lfq{\leavevmode\raise.3ex\hbox{$\scriptscriptstyle
\langle\!\langle$}\thinspace}
\def\rfq{\leavevmode\thinspace\raise.3ex\hbox{$\scriptscriptstyle
\rangle\!\rangle$}}

\auteurcourant{DOMINIQUE FOATA {\sevenrm AND} GUO-NIU
HAN}
\titrecourant{TANGENT AND SECANT $q$-DERIVATIVE
POLYNOMIALS}

\catcode`\@=12


\def\ides{\mathop{\rm ides}\nolimits}

\def\imaj{\mathop{\rm imaj}\nolimits}
\def\inv{\mathop{\rm inv}\nolimits}
\def\Ligne{\mathop{\rm Ligne}\nolimits}
\def\Iligne{\mathop{\rm Iligne}\nolimits}

\def\Sin{\mathop{\rm Sin}\nolimits}
\def\Cos{\mathop{\rm Cos}\nolimits}
\def\Sec{\mathop{\rm Sec}\nolimits}
\def\Tan{\mathop{\rm Tan}\nolimits}

\def\min{\mathop{\rm min}\nolimits}

\def\starAS{\mathop{{}^*\kern-3pt{\rm AS}}\nolimits}
\def\ASprimstar#1{\mathop{{\rm
AS}^{\lower3pt\hbox{$'$}*}}_{\kern-7pt#1}\nolimits}

\def\Dun{\displaystyle
\mathop{D}_{\scriptscriptstyle<}\kern-2pt
{}_q}
\def\Ddeux{\displaystyle
\mathop{D}_{\scriptscriptstyle\ge}\kern-2pt {}_q}

\def\qed{\quad\raise -2pt\hbox{\vrule\vbox to 10pt{\hrule width 4pt
\vfill\hrule}\vrule}}

\def\article#1|#2|#3|#4|#5|#6|#7|
    {{\leftskip=7mm\noindent
     \hangindent=2mm\hangafter=1
     \llap{{\tt [#1]}\hskip.35em}{#2}.\quad
     #3, {\sl #4}, {\bf #5} ({\oldstyle #6}),
     pp.\nobreak\ #7.\par}}
\def\livre#1|#2|#3|#4|
    {{\leftskip=7mm\noindent
    \hangindent=2mm\hangafter=1
    \llap{{\tt [#1]}\hskip.35em}{#2}.\quad
    {\sl #3}, #4.\par}}
\def\divers#1|#2|#3|
    {{\leftskip=7mm\noindent
    \hangindent=2mm\hangafter=1
     \llap{{\tt [#1]}\hskip.35em}{#2}.\quad
     #3.\par}}

\noindent

June 18, 2012  
\bigskip\bigskip
\centerline{\bf Multivariable Tangent and Secant} 
\smallskip
\centerline{\bf $q$-derivative Polynomials}

\bigskip
\centerline{\sl Dominique Foata, Guo-Niu Han}
\footnote{}{
}

\bigskip\bigskip

{
\eightpoint
\narrower
\noindent
{\bf Abstract}.\quad
The derivative polynomials introduced by Knuth and
Buckholtz in their calculations of the tangent and secant
numbers are extended to a multivariable
$q$--environment. The $n$-th $q$-derivatives of the
classical $q$-tangent and 
$q$-secant are each given
two polynomial expressions. 
The first polynomial expression is 
indexed by triples of integers,
the second by compositions of integers. 
The
functional relation between those two classes is fully given
by means of combinatorial techniques. Moreover, those
polynomials are proved to be generating functions for
so-called $t$-permutations by 
multivariable statistics. By giving special values to
those polynomials we recover classical $q$-polynomials
such as the Carlitz $q$-Eulerian polynomials and the
$(t,q)$-tangent and -secant analogs recently introduced.
They also provide 
$q$-analogs for the Springer numbers.
Finally, the $t$-compositions used in this paper furnish a combinatorial interpretation to one of the Fibonacci triangles.

}

\bigskip
\centerline{\bf Summary}

{\eightpoint
\medskip
\halign{\indent\indent\hfil #.&\ #\hfil\cr
1&Introduction\cr
\omit&1. The derivative polynomials\cr
\omit&2. Towards a multivariable $q$-analog\cr
\omit&3. The numerical and combinatorial background\cr
\omit&4. The underlying statistics\cr
\omit&5. The main results\cr
2&A detour to the theory of $q$-trigonometic functions\cr
3&Transformations on $t$-permutations\cr
4&Proof of Theorem 1.1\cr
5&Proof of Theorem 1.2\cr
6&Proof of Theorem 1.3\cr
7&More on $q$-trigonometric functions\cr
8&Proof of Theorem 1.4 \cr
9&Proof of Theorem 1.5 \cr
10&Specializations\cr
\omit&1. The first column of $(A_{n,k,a,b}(q))$,
$(B_{n,k,a,b}(q))$\cr
\omit&2. The super-diagonal $(n,n+1)$ of the matrix
$(A_{n,k,a,b}(q))$\cr
\omit&3. The subdiagonal $(n,n-1)$ of the matrix
$(A_{n,k,a,b}(q))$\cr
\omit&4. The subdiagonal $(n,n-2)$ of the matrix
$(B_{n,k,a,b}(q))$\cr
\omit&5. Two $q$-analogs of the Springer numbers\cr
\omit&6. $t$-compositions and Fibonacci triangle\cr
\omit&7. Further comment\cr
11&Tables\cr
\omit&References\cr
}

}

\vfill\eject
\bigskip\bigskip
\centerline{\bf 1. Introduction}  

\medskip
Back in 1967, Knuth and Buckholtz [KB67] devised a clever
method for computing the tangent and secant numbers
$T_{2n+1}$ and $E_{2n}$ for large values of the subscripts~$n$.
For that purpose they introduced two sequences of polynomials,  
referred to as {\it derivative polynomials}. A few years
later, Hoffman [Ho95, Ho99] calculated the exponential generating
functions for those polynomials, and found a combinatorial
interpretation for their coefficients, in terms of 
so-called {\it snakes}. The goal of this paper is to obtain
a {\it multivariable $q$-analog} of all those results. We 
first recall the contributions made by those authors, then,
introduce the $q$-environment that makes it possible to derive
a handy algebra for these new $q$-derivative polynomials.

\medskip
1.1. {\it The derivative polynomials}.\quad
Recall that {\it tangent} and {\it secant numbers} $T_{2n+1}$ 
and $E_{2n}$ occur as coefficients in the
Taylor expansions of $\tan u$ and $\sec u$: 
$$
\leqalignno{\quad\tan u&=\sum_{n\ge 0} {u^{2n+1}\over
(2n+1)!}T_{2n+1}&(1.1)\cr
&={u\over 1!}1+{u^3\over 3!}2+{u^5\over 5!}16+{u^7\over
7!}272+ {u^9\over 9!}7936+\cdots\cr  
\qquad\sec u={1\over \cos u}&=\sum_{n\ge 0}
{u^{2n}\over (2n)!}E_{2n}&(1.2)\cr
&=1+{u^2\over 2!}1+{u^4\over 4!}5+{u^6\over 6!}61+{u^8\over
8!}1385+ {u^{10}\over 10!}50521+\cdots\cr  }
$$
See, e.g.,
[Ni23, p.~177-178], [Co74, p.~258-259]. 

Let $(A_n(x))$ $(n\ge 1)$ be the sequence of polynomials
defined by
$$
A_0(x)=x,\quad A_{n+1}=(1+x^2)\,DA_n(x),
$$
with $D$ being the differential operator. When writing
$A_n(x)=\smash{\sum\limits_{m\ge 0}}a(n,m)x^m$, the
coefficients
$(a(n,m))$ satisfy the recurrence
$$
a(0,m)=\delta_{1,m},\quad
a(n+1,m)=(m-1)a(n,m-1)+(m+1)a(n,m+1).\leqno(1.3)
$$
The $a(n,m)$'s form a triangle of integral
numbers (see Table~1 in Section 11, where the numbers
$a(n,m)$ are reproduced in boldface), now  registered as
the sequence A101343 in Sloane's Encyclopedia of Integer
Sequences [Sl06] with an abundant bibliography.
Knuth and Buckholtz [KB67] showed that the 
$n$-th derivative $D^n \tan u$ was equal to the
polynomial
$$
D^n \tan u=\sum_{m\ge 0}a(n,m)\tan^m u.\leqno(1.4)
$$

The same two authors also introduced the sequence
$(b(n,m))$ by 
$$
b_{0,m}=\delta_{0,m},\quad
b(n+1,m)=mb(n,m-1)+(m+1)b(n,m+1).
$$

\goodbreak\noindent
Again, Knuth and Buckholtz [KB67] showed that the $n$-th
derivative of $\sec u$ could be expressed as
$$
D^n\sec u:=\sum_{m\ge 0}b(n,m)\tan^m u\,\sec u.
\leqno(1.5)
$$

The triangle of numbers $(b(n,m))$ also appears in
Sloane's Enclycopedia [Sl06] with reference A008294. The
first values of the numbers $b(n,m)$ are reproduced in
Table~1 in plain type (not bold). From their very
definitions the $a(n,m)$'s and $b(n,m)$'s can
be imbricated in the same table, as done in Table~1. The
meanings of the entries to the right of the table will  be
further explained. 

The exponential generating functions for the polynomials
$$
A_n(x):=\sum_{m\ge 0}a(n,m)x^m\quad{\rm and}\quad
B_n(x):=\sum_{m\ge 0}b(n,m)x^m,
$$ 
called {\it derivative polynomials},
have been derived by
Hoffman [Ho95] in the form
$$\leqalignno{
\sum_{n\ge 0}A_n(x){u^n\over n!}
&={x+\tan u\over 1-x\,\tan u};&(1.6)\cr
\sum_{n\ge 0}B_n(x){u^n\over n!}
&={1\over \cos u-x\sin u}.&(1.7)\cr}
$$
Those two exponential generating functions
and recurrences for the $a(n,m)$'s and $b(n,m)$'s have also
been obtained by other people in different contexts, in particular,
by Carlitz and Scoville [CS72], Fran\c con [Fr78].

By plugging $x=1$ in (1.6) the right-hand side becomes
$\tan 2u+\sec 2u$, so that the sum $\sum_m a(n,m)$
is equal to $2^nE_n$, if $n$ is even, and to~$2^nT_n$ if
$n$ is odd. Likewise, (1.7) yields
$\sum_nB_n(1)u^n/n!=1/(\cos u-\sin u)$, which is the
exponential generating function for the so-called {\it
Springer numbers} (1, 1, 3, 11, 57, 361, 2763, \dots)
(see [Sp71], [Du95]) originally considered by Glaisher [Gl98, Gl99,
Gl14], as noted in Sloane's Encyclopedia, under reference
A001586. In Table~1 we have indicated the values of the
row sums $\sum_ma(n,m)$ and $\sum_mb(n,m)$ to the right.

Finally, the combinatorial interpretations of the $a(n,m)$'s
and $b(n,m)$'s are due to Hoffman in a later paper
[Ho99]. An equivalent interpretation is also due to
Josuat-Verg\`es [Jo11]. Both authors use the word {\it
snake} of length~$n$, a notion made popular by Arnold
[Ar92, Ar92a] in the study of morsification of singularities,
to designate each word
$w=x_1x_2\cdots x_n$, whose letters are integers,
positive or negative, with the further property that
$x_1>x_2$, $x_2<x_3$, $x_3>x_4$, \dots\ in an
alternating way and
$|x_1|\,|x_2|\,\cdots\,|x_n|$ is a permutation of
$1\,2\,\cdots\,n$. Note that Josuat-Verg\`es, Novelli 
and Thibon [JNT12] have recently developed an algebraic 
combinatorics of snakes from the Hopf algebra
point of view. In [Ho99] and [Jo11] Hoffman and Josuat-Verg\`es show
that $A_n(x)$ is the generating polynomial for the set of all
snakes of length~$n$ by the number of sign changes; they
also prove an analogous result for $B_n(x)$.

\goodbreak
\medskip
1.2. {\it Towards a multivariable $q$-analog}.\quad
In parallel with (1.4) and (1.5) the $q$-{\it derivative
operator} $D_q$ (see [GR90, p.~22]), as well as the
$q$-{\it analogs} of
tangent and secant (see [St76], [AG78], [AF80], [Fo81], [St97, p.
148-149], [St10]) are to be introduced.  The first problem is to see
whether the $n$-th $q$-derivatives of those $q$-analogs
can be expressed as polynomials in those functions, and
they can!  But, contrary to formulas (1.4) and (1.5), those
$n$-th $q$-derivatives have {\it several} polynomial
forms. As will be seen, {\it two} such polynomial forms are
derived in this paper for each $q$-analog of tangent and
secant. The second problem is to work out an appropriate
algebra for the polynomials involved that must appear as
natural {\it multivariable}
$q$-analogs of the entries $a(n,m)$ and
$b(n,m)$. 

Before stating the first results of this paper we
recall a few basic notions on $q$-Calculus.
The $t$-{\it ascending factorial} in a variable~$t$ is
traditionally defined by
$$\leqalignno{
(t;q)_n&:=\cases{1,&if $n=0$;\cr
(1-t)(1-tq)\cdots (1-tq^{n-1}),&if $n\ge 1$;\cr}\cr
\noalign{\hbox{in its finite version and}}
\noalign{\vskip-2pt}
(t;q)_\infty&:=\textstyle \lim_n(t;q)_n=\prod\limits_{n\ge
0}(1-tq^n); \cr
\noalign{\vskip-4pt}
}
$$
in its infinite version. By $q$-series it is meant a series
of the form $f(u)=\sum_{n\ge 0}f(n;q) u^n/(q;q)_n$,
whose coefficients $f(n;q)$ belong to some 
ring [AAR00, chap.~10]. The $q$-{\it derivative
operator} for fixed~$q$ used below is defined by
$$
D_qf(u):={f(u)-f(qu)\over u},
$$
instead of the traditional $(f(u)-f(qu))/(u(1-q))$.

At this stage we just have to note (see Section~2) that
there is {\it only one $q$-tangent} attached to the
classical Jackson definitions [Ja04] of the $q$-sine and
$q$-cosine, namely,
$$\leqalignno{\noalign{\vskip-10pt}
\tan_q(u)&:={\sum\limits_{n\ge
0}(-1)^nu ^{2n+1}/(q;q)_{2n+1}
\over \sum\limits_{n\ge 0}(-1)^n u^{2n}/(q;q)_{2n}},
\cr
}
$$
but, as it has been rarely noticed, 
there are {\it two $q$-secants}
$$\leqalignno{
\sec_q(u)&:
={1\over\sum\limits_{n\ge 0} (-1)^n u^{2n} /
(q;q)_{2n}};\cr
\Sec_q(u)&:
={1\over\sum\limits_{n\ge 0} (-1)^n  q^{n(2n-1)}u^{2n} /
(q;q)_{2n}}.\cr}
$$

1.3. {\it The numerical and combinatorial background}.\quad
Following Andrews [An76, chap. 4] a {\it composition} of a
positive integer~$n$ is defined to be a sequence 
${\bf c}=(c_0,c_1,\ldots, c_m)$ of {\it nonnegative}
integers such that $c_0+c_1+\cdots +c_m=n$, 
with the restriction that only $c_0$ and $c_m$ can be zero; the
$c_i$'s are the {\it parts} of the composition and $m+1$
is the {\it number of parts}, denoted by $\mu\,{\bf
c}+1$.
A composition ${\bf c}=(c_0,c_1,\ldots, c_m)$ of~$n$ is
said to be a $t$-{\it composition}, if the following two
conditions hold:

\goodbreak
(1) either $m=0$, so that  ${\bf c}=(c_0)$ and $c_0=n$ is
an {\it odd} integer, or $m\ge 1$ and both $c_0$, $c_m$
are {\it even};

(2) if $m\ge 2$, then all the parts $c_1$, $c_2$, \dots~,
$c_{m-1}$ are {\it odd}.

\noindent For each $n\ge 1$ the set of all
$t$-compositions of~$n$ is denoted by $\Theta_n$. It is
further assumed that
$\Theta_0$ consists of a unique empty composition
denoted by~$(0,0)$. 

\smallskip
The first $t$-compositions are the following:

\noindent
$\Theta_0$: \quad $(0,0)$;

\noindent
$\Theta_1$: \quad $(1)$,\quad $(0,1,0)$;

\noindent
$\Theta_2$: \quad $(0,2)$,\quad $(2,0)$,\quad $(0,1,1,0)$;

\noindent
$\Theta_3$: \quad $(3)$,\quad $(0,1,2)$,\quad
$(2,1,0)$,\quad $(0,3,0)$,\quad $(0,1,1,1,0)$;

\noindent
$\Theta_4$: \quad $(0,4)$,\ $(2,2)$, $(4,0)$,\ 
$(0,3,1,0)$,\  $(0,1,3,0)$,\  $(0,1,1,2)$,\ 
$(2,1,1,0)$,\hfil\break
\null\kern28pt  $(0,1,1,1,1,0)$.
\smallskip
\goodbreak

A $t$-composition ${\bf c}=(c_0, c_1, \ldots, c_{m-1}, c_m)$ from $\Theta_n$ such that $c_m=0$
is called an {\it $s$-composition}. We write ${\bf c}^-=(c_0, c_1, \ldots, c_{m-1})$.
The subset of $\Theta_n$ of all $s$-compositions is denoted by $\Theta_n^-$. 
\medskip

For each composition 
${\bf c}=(c_0,c_1,\ldots,c_m)$ of
$n\ge 1$ let $\tan_q(q^{\bf c}u):=1$ if $m=0$ and
for $m\ge 1$, whatever $n\ge 1$, let
$$
\tan_q(q^{\bf c}u):=
\tan_q(q^{c_0}u)\tan_q(q^{c_0+c_1}u)\cdots
\tan_q(q^{c_0+c_1+\cdots+c_{m-1}}u).
$$
Note that the rightmost
part~$c_m$ does not occur in the previous expression.
Also let $\rho\,{\bf c}:=(c_m,\ldots,c_1,c_0)$ denote the mirror-image of~${\bf c}$. 

\medskip

A word $w=y_1y_2\cdots y_m$, whose letters are positive integers,
is said to be {\it falling
alternating}, or simply {\it alternating} (resp. {\it rising
alternating}), if $y_1>y_2$, $y_2<y_3$, $y_3>y_4$,
\dots\  (resp. if $y_1<y_2$, $y_2>y_3$, $y_3<y_4$,
\dots) in an alternating manner. The notion goes back to
D\'esir\'e Andr\'e [An79, An81], who showed that the
number of falling (resp. rising) alternating permutations
$\sigma=\sigma(1)\sigma(2) \cdots \sigma(n)$ of
$1\,2\,\ldots\,n$ is equal to $T_n$ when $n$ is odd, and
to $E_n$ when~$n$ is even.
The length of each word~$w$ is denoted by $\lambda w$ and
the {\it empty word} is the word of length~0, denoted
by~$\epsilon$. 

\medskip
{\it Definition}.\quad 
A $t$-{\it permutation} of order~$n$ is defined to be a
nonempty sequence
$w=(w_0,w_1,\ldots,w_m)$ of words having the
properties:

(i) the juxtaposition product $w_0w_1\cdots w_m$ is  a
permutation of $12\cdots n$;

(ii) either $m=0$ and $w_0$ is rising alternating of
{\it odd} length, or $m\ge 1$ and then~$w_0$
is rising alternating of {\it even} length and~$w_m$ is
(falling) alternating of {\it even} length;

\goodbreak
(iii) if $m\ge 2$, then all the components $w_1$, $w_2$,
\dots~, $w_{m-1}$ are (falling) alternating of {\it odd}
length.

\medskip
For each $n\ge 0$ the set of
all $t$-permutations  of order~$n$ will
be denoted by~${\cal T}_n$. The first $t$-permutations are the following:
\halign{#\hfil\quad&#\hfil\cr
${\cal T}_0$:&$(\epsilon,\epsilon)$
;\cr 
${\cal T}_1$:&$(1)$;\quad $(\epsilon,1,\epsilon)$;\cr
${\cal T}_2$:&$(\epsilon,21),\ (12,\epsilon)$;
\quad $(\epsilon,2,1,\epsilon)$,\ 
$(\epsilon,1,2,\epsilon)$;\cr}

\goodbreak
\halign{#\hfil\quad&#\hfil\cr
${\cal T}_3$:&$(132),\  (231)$;\cr 
&$(\epsilon,3,21),\ (\epsilon,312,\epsilon),\
(\epsilon,213,\epsilon),\  (12,3,\epsilon),
(\epsilon,1,32),\
(23,1,\epsilon)$,\cr
&$(13,2,\epsilon),\ (\epsilon,2,31)$;\cr 
&$(\epsilon,1,2,3,\epsilon),\
(\epsilon,1,3,2,\epsilon),\  (\epsilon,2,1,3,\epsilon),\ 
(\epsilon,2,3,1,\epsilon)$,\cr 
&$(\epsilon,3,1,2,\epsilon),\ (\epsilon,3,2,1,\epsilon)$.\cr}

\medskip
1.4. {\it The underlying statistics}.\quad 
The number of factors in each $t$-permutation 
$w=(w_0,w_1,\ldots,w_m)$, minus one,
is denoted by $\mu\, w=m$. If $w$ is of order~$n$, the
sequence  $(\lambda w_0,\lambda w_1,\ldots,\lambda w_m)$
is a $t$-composition of~$n$, denoted by $\Lambda
w$. For each $n\ge 0$ the set
of all $t$-permutations $w$ from ${\cal T}_n$, such that $\mu\, w=m$,
with $0\leq m \leq n+1$ and $m\equiv n+1 \bmod 2$, will be denoted 
by ${\cal T}_{n,m}$. 

Recall that the {\it ligne of route} of a
permutation $\sigma=\sigma(1)\sigma(2)\cdots\sigma(n)$
of $12\cdots n$, denoted by $\Ligne\sigma$,  is
defined to be the set of all~$i$ such that $1\le i\le n-1$
and
$\sigma(i)>\sigma(i+1)$; also, the {\it inverse ligne of
route}, $\Iligne\sigma$, to be the set of all $\sigma(i)$
such that $\sigma(j)=\sigma(i)+1$ for some $j\le i-1$.
In an equivalent manner, $\Iligne\sigma=\Ligne \sigma^{-1}$.
Next, let
$$
\leqalignno{\ides\sigma&=\#\Iligne\sigma;\cr
\imaj\sigma&=\sum_i\sigma(i)\quad (\sigma(i)\in
\Iligne\sigma);\cr}
$$
and let $\inv\sigma$ be the {\it number of inversions}
of~$\sigma$, as being the number of pairs $(i,j)$ such
that $1\le i<j\le n$ and $\sigma(i)>\sigma(j)$. 

The {\it inverse ligne of route}, $\Iligne w$, of a
$t$-permutation   $w=(w_0,w_1,\ldots,w_m)$ is
then defined by
$$
\leqalignno{
\Iligne w&:=\Iligne (w_0w_1\cdots
w_m);\cr
\noalign{\hbox{and the {\it number of inversions},
$\inv w$, by}}
\inv w&:=\inv (w_0w_1\cdots w_m).\cr
}
$$
This makes sense, as the latter juxtaposition product is a
permutation, say,~$\sigma$ of $12\cdots n$ if $w\in
{\cal T}_n$. Finally, let $\min w:=a$ if 1 is
a letter in~$w_a$.

\medskip
For example, with the $t$-permutation
$w=({4}\,
{5},\,11\,{1}\,{\bf3},\, {\bf 1\!0}\, {7}\,
{\bf9},\,{\bf6},\, {\bf8}\, {\bf2})$, the elements of $\Iligne
w$ being reproduced in boldface, we have:
$\ides w=6$, $\imaj w=\bf 3+10+9+6+8+2=38$,
$\inv w=27$ and
$\min w=1$.

\medskip
1.5. {\it The main results}.\quad
For each triple $(k,a,b)$ let ${\cal T}_{n,k,a,b}$ 
denote the set of all $t$-permutations $w$ from
${\cal T}_n$ such that
$\ides w=k$, $\min w=a$ and $a+b=\mu\,w$.

\proclaim Theorem 1.1 {\rm (Multivariable $q$-analog of (1.4))}. 
Let 
$$
A_{n,k,a,b}(q)=\sum_{w\in {\cal T}_{n,k,a,b}}
q^{\imaj w}
\leqno(1.8)
$$
be the generating polynomial for the set
${\cal T}_{n,k,a,b}$ by the statistic ``\/$\imaj$." 
Then
$$
D_q^n\tan_q(u)=\sum_{k,a,b}A_{n,k,a,b}(q)\,
(\tan_q(q^{k+1}u))^b\,(\tan_q(q^ku))^a;\leqno(1.9)
$$
where  $0\le k\le n-1$ and $0\le a+b\le n+1$.

\goodbreak
In the same manner, for each $(k,a,b)$ let
${\cal T}_{n,k,a,b+1}^-$ denote the set of all 
$t$-permutations $w=(w_0, w_1, \ldots, w_m, w_{m+1})$ in ${\cal T}_{n,k,a,b+1}$ 
such that $w_{m+1}=\epsilon$. 

\proclaim Theorem 1.2 {\rm (Multivariable $q$-analog of (1.5))}. Let 
$$
B_{n,k,a,b}(q)=\sum_{w\in {\cal T}_{n,k,a,b+1}^-}
q^{\imaj w}
\leqno(1.10)
$$
be the generating polynomial for the set
${\cal T}_{n,k,a,b+1}^-$ by ``\/$\imaj$." Then
$$\leqalignno{
\qquad D_q^n\sec_q(u)&=\sum_{k,a,b}B_{n,k,a,b}(q)
\bigl(\tan_q(q^{k+1}u)\bigr)^b\,
\sec_q(q^{k+1}u)\,\bigl(\tan_q(q^{k}u)\bigr)^a;&(1.11)\cr
\qquad D_q^n\Sec_q(u)&=\sum_{k,a,b}q^{n(n-1)/2}B_{n,n-1-k,b,a}(q^{-1}) &(1.12)\cr
\noalign{\kern -10pt}
&\qquad\qquad \times \bigl(\tan_q(q^{k+1}u)\bigr)^b\,
\Sec_q(q^{k}u)\,\bigl(\tan_q(q^{k}u)\bigr)^a.\cr}
$$
where $0\le k\le n-1$ and $0\le a+b\le n$.

\proclaim Theorem 1.3 {\rm (Composition $q$-analogs of (1.4) and (1.5))}. 
For each $n\ge 1$ and each
$t$-composition ${\bf c}$ of $n$ let $A_{n,{\bf c}}(q)$ be the polynomial
$$\leqalignno{
A_{n,{\bf c}}(q)&=\sum_{w\in {\cal T}_n,\Lambda w={\bf c}}
q^{\inv w}.&(1.13)\cr
\noalign{\hbox{Then,}}
D_q^n \tan_q(u) &=\sum_{{\bf c}\in \Theta_n} A_{n,{\bf c}}(q)\,\tan_q(q^{\bf c}u);&(1.14)\cr
D_q^n\sec_q(u)&=\sum_{{\bf c}\in\Theta_n^{-}} 
  A_{n,{\bf c}}(q)\,\tan_q(q^{{\bf c}^{-}}u)\,\sec_q(q^n
u);&(1.15)\cr
D_q^n\Sec_q(u)&=\sum_{{\bf c}\in\Theta_n^{-}} q^{n(n-1)/2} 
	A_{n,{\bf c}}(q^{-1})\,\tan_q(q^{\rho {\bf c}^{-}}u)\,
\Sec_q(u).&(1.16)\cr}
$$

\medskip
{\it Remark}.\quad
The polynomials $A_{n,k,a,b}(q)$ are not uniquely defined by identity (1.9). On the contrary, both (1.11) and (1.12) uniquely define the polynomials $B_{n,k,a,b}(q)$ and $A_{n,{\bf c}}(q)$.

\medskip 
Polynomials indexed by triples $(k,a,b)$ and those by compositions~$\bf c$ are related to each other, as indicated in the next theorem.

\proclaim Theorem 1.4.
We have
$$
\leqalignno{\noalign{\vskip-4pt}
\sum_{k\ge
0,\,a+b=m}A_{n,k,a,b}(q) &=\sum_{{\bf c}\in\Theta_{n,m}}A_{n,{\bf
c}}(q);&(1.17)\cr
\noalign{\vskip-2pt}
\sum_{ k\ge
0,\,a+b=m}B_{n,k,a,b}(q)&=\sum_{{\bf c}\in\Theta_{n,m+1}^-}A_{n,{\bf
c}}(q).&(1.18)\cr
\noalign{\vskip-2pt}}
$$

Now, form the generating polynomials 
$$
\leqalignno{
A_n(x,q)&:=\sum_{m\ge 0}x^m\kern-6pt
\sum_{ k\ge 0,\,a+b=m}\kern-12pt
A_{n,k,a,b}(q)
=\sum_{m\ge 0}x^m\sum_{{\bf c}\in\Theta_{n,m}}A_{n,{\bf
c}}(q);\cr
B_n(x,q)&:=\sum_{m\ge 0}x^m\kern-6pt
\sum_{ k\ge
0,\,a+b=m}\kern-12ptB_{n,k,a,b}(q)
=\sum_{m\ge 0}x^m\sum_{{\bf c}\in\Theta_{n,m+1}^-}A_{n,{\bf
c}}(q).\cr}$$

\proclaim Theorem 1.5. The factorial generating functions
for the polynomials $A_n(x,q)$ and $B_n(x,q)$ are given by:
$$\leqalignno{
\qquad\sum_{n\ge 0}A_n(x,q){u^n\over (q;q)_n}
&=\tan_q(u)+\sec_q(u)(1-x\tan_q(u))^{-1}x\Sec_q(u);&(1.19)\cr
\sum_{n\ge 0}B_n(x,q){u^n\over (q;q)_n}
&=\sec_q(u)(1-x\tan_q(u))^{-1}.&(1.20)\cr
}
$$

Those two formulas, derived in Section~9, are true
$q$-analogs of Hoffman's identities (1.6) and (1.7), as
the latter ones can be rewritten as: 
$$
\leqalignno{\noalign{\vskip-7pt}
\qquad\sum_{n\ge 0}A_n(x){u^n\over n!}
&=\tan(u)+\sec(u)(1-x\tan(u))^{-1}x\sec(u),\cr
\sum_{n\ge 0}B_n(x){u^n\over n!}
&=\sec(u)(1-x\tan(u))^{-1}.&(1.21)\cr
}
$$

The four-variable polynomials
$A_n(t,x,y,q)=\sum\limits_{k,a,b}
A_{n,k,a,b}(q)t^kx^ay^b$
and
$B_n(t,x,y,q):=\sum\limits_{k,a,b}B_{n,k,a,b}(q)t^kx^ay^b$,
that may be considered as multivariable $q$-analogs of the
entries $a(n,m)$ and $b(n,m)$, have several
interesting {\it specializations} studied in Section~10.
First,
$tA_n(t,0,0,q)$ and $tB_n(t,0,0,q)$ are shown to be the
$(t,q)$-analogs $T_n(t,q)$ and $E_n(t,q)$ of tangent and
secant numbers, only defined so far [FH11] by their factorial
generating functions
$\sum_{n}T_{2n+1}(t,q)u^{2n+1}/
(t;q)_{2n+2}$ and
$\sum_{n}E_{2n}(t,q)u^{2n}/ (t;q)_{2n+1}$. The
recurrences of the polynomials
$A_{n,k,a,b}(q)$ and $B_{n,k,a,b}(q)$ provide a handy
method for calculating them. We also prove that
the polynomial
$A_n(t,x,q):=\sum_{k,a}A_{n,k,a,n+1-a}(q)t^kx^a$ is a {\it
refinement} of the  Carlitz
$q$-analog [Ca54] $A_n(t,q)$ of the {\it Eulerian polynomial}, 
	as $A_n(t,1,q)=A_n(t,q)$, with an explicit
combinatorial interpretation. Referring to 
Tables~2 and~3 of the polynomials $A_{n,k,a,b}(q)$ and 
$B_{n,k,a,b}(q)$ displayed at the end of the paper, it is
shown that the sum of the polynomials occurring in each
box along the two top diagonals can be
explicitly evaluated. 
Furthermore, the Springer numbers
are given two $q$-analogs. Finally, generating functions and
recurrence relations are provided for both $t$- and
$s$-compositions.

\bigskip

\centerline{\bf 2. A detour to the theory of 
$q$-trigonometric functions}

\medskip
By means of the $q$-binomial
theorem  ([GR90, \S\kern2pt 1.3]; [AAR00, \S\kern2pt 10.2])
we can express the {\it first $q$-exponential}
$e_q(u)$ and the {\it second $q$-exponential} $E_q(u)$,
either as an infinite $q$-series, or an infinite product:
$$\displaylines{\rlap{(2.1)}\hfill
e_q(u)=
\sum_{n\ge 0} {u^n\over
(q;q)_n}={1\over (u;q)_\infty};\hfill\cr
\rlap{(2.2)}\hfill
E_q(u)=\sum_{n\ge 0}q^{n(n-1)/2}
{u^n\over (q;q)_n}=(-u;q)_\infty,\hfill\cr}
$$
two results that go back to Euler [Eu48].
As already done by Jackson [Ja04] (also see [GR90, p.~23]),
they both serve to define the $q$-trigonometric functions
$q$-sine and $q$-cosine:
$$\displaylines{
\sin_q(u):={e_q(iu)-e_q(-iu)\over 2i}=
\sum_{n\ge 0} (-1)^n {u^{2n+1} \over (q;q)_{2n+1}}\,;\cr
\cos_q(u):={e_q(iu)+e_q(-iu)\over 2}=
\sum_{n\ge 0} (-1)^n {u^{2n} \over (q;q)_{2n}}\,;\cr
\Sin_q(u):={E_q(iu)-E_q(-iu)\over 2i}\qquad{\rm
and}\qquad
\Cos_q(u):={E_q(iu)+E_q(-iu)\over 2}.\cr
}
$$
We can also define:
$\displaystyle
\tan_q(u)\!:=\!{\sin_q(u)/\cos_q(u)}$  and
$\displaystyle\Tan_q(u)\!:=\!{\Sin_q(u)/ \Cos_q(u)}$; but,
as $\sin_q(u)\Cos_q(u)-\Sin_q(u)\cos_q(u)=0$, we have:
$$\displaylines{
\qquad
\tan_q(u)={\sin_q(u)\over \cos_q(u)}
={\Sin_q(u)\over \Cos_q(u)}=\Tan_q(u),\qquad\cr}
$$
so that there is {\it only one} $q$-tangent. However, there
are two $q$-secants:
$$\leqalignno{
\sec_q(u)&:={1\over\cos_q(u)}
={1\over\sum\limits_{n\ge 0} (-1)^n u^{2n} /
(q;q)_{2n}};\cr
\Sec_q(u)&:={1\over\Cos_q(u)}
={1\over\sum\limits_{n\ge 0} (-1)^n  q^{n(2n-1)}u^{2n} /
(q;q)_{2n}}.\cr}
$$

\proclaim Theorem 2.1. The $q$-derivatives of the series
$\tan_q(u)$, $\sec_q(u)$, $\Sec_q(u)$ can be evaluated as
follows:
$$\leqalignno{
D_q\tan_q(u)&=1+\tan_q(u)\tan_q(qu);&(2.3)\cr
D_q\sec_q(u)&=\sec_q(qu)\,\tan_q(u);&(2.4)\cr
D_q\Sec_q(u)&=
\Sec_q(u)\tan_q(qu).&(2.5)\cr}
$$

{\it Proof}.\quad These three identities can be proved,
either by working with $e_q(u)$ and $E_q(u)$, when
expressed as infinite $q$-series, or by using the infinite
products appearing in (2.1) and (2.2). We choose the former
way, because of its closeness to the traditional
trigonometric calculus.

First, 
$$
\leqalignno{e_q(\alpha u)-e_q(\alpha qu)
&=\sum_{n\ge 1}{(\alpha
u)^n(1-q^n)\over (q;q)_n}
=\alpha u\sum_{n\ge 1}
{(\alpha u)^{n-1}\over (q;q)_{n-1}}=\alpha u\, e_q(\alpha
u);\cr
E_q(\alpha u)-E_q(\alpha qu)&=\sum_{n\ge
1}{(\alpha u)^nq^{n(n-1)/2}(1-q^n)\over (q;q)_n}\cr
&\kern1.5cm{}=\alpha
u\sum_{n\ge 1}{(\alpha
uq)^{n-1}q^{(n-1)(n-2)/2}\over (q;q)_{n-1}}=\alpha u\,
E_q(\alpha qu).\cr
\noalign{\hbox{Hence,}}
D_q\,e_q(\alpha u)&=\alpha\,e_q(\alpha u);
\qquad D_q\,E_q(\alpha u)=\alpha \,E_q(\alpha
qu).\cr }$$

\goodbreak\noindent
Next, by applying $D_q$ to the familiar identities:
$e_q(iu)=\cos_q(u)+i\sin_q(u)$ and
$E_q(iu)=\Cos_q(u)+i\Sin_q(u)$, we get
$$\leqalignno{
D_q\cos_q(u)&=-\sin_q(u);\quad 
D_q\sin_q(u)=\cos_q(u);\cr
D_q\Cos_q(u)&=-\Sin_q(u);\quad 
D_q\Sin_q(u)=\Cos_q(u).\cr}
$$
Finally, we take advantage of the next formula that yields
an expression for the $q$-derivative of a ratio $f(u)/g(u)$
of two $q$-series
$$\leqalignno{
D_q\,{f(u)\over g(u)}
&={g(qu)\,D_qf(u)-f(qu)\,D_qg(u)
\over g(u)g(qu)},\cr
\noalign{\hbox{and use it for the ratios
$\sin_q(u)/\cos_q(u)$, $1/\cos_q(u)$, $1/\Cos_q(u)$ to
get:}} 
D_q\tan_q(u)
&={\cos_q(qu)\cos_q(u)-\sin_q(qu)(-\sin_q(u))\over
\cos_q(u)\cos_q(qu)} =1+\tan_q(u)\tan_q(qu);\cr
D_q\sec_q(u)&={\sin_q(u)\over \cos_q(u)\cos_q(qu)}=
\sec_q(qu)\tan_q(u);\cr
D_q\Sec_q(u)&={\Sin_q(qu)\over \Cos_q(u)\Cos_q(qu)}=
\Sec_q(u)\tan_q(qu).\qed\cr
}
$$

\medskip
Identities (2.4) and (2.5) above show a certain duality between the $q$-derivatives of $\sec u$ and $\Sec u$, which is to be explored by using the base~$q^{-1}$ instead of~$q$. 
First, for each $q$-series~$f(u)$ let $Qf(u):=f(qu)$ and $Uf(u):=(1/u)f(u)$. 
We have the identities:
$$
UQ=q\,QU,\qquad Q^{-1}U=q\,UQ^{-1}.
$$
Next, the $q$-difference operators $D_q$ and
$D_{q^{-1}}$, the latter being defined by $D_{q^{-1}}f(u):=(1/u)(f(u)-f(q^{-1}u))$, also read
$$\displaylines{\rlap{(2.6)}\hfill
D_{q}=U(I-Q);\qquad
D_{q^{-1}}=U(I-Q^{-1}).\hfill\cr
\noalign{\hbox{Hence,}}
\rlap{(2.7)}\hfill D_q=-D_{q^{-1}}Q.\hfill\cr}
$$
 Next, we have the relations:
$$\displaylines{\rlap{(2.8)}\hfill
e_{q^{-1}}(u)=QE_q(-u),\qquad
E_{q^{-1}}(u)=Qe_q(-u).\hfill\cr
\rlap{(2.9)}\hfill
\sin_{q^{-1}}(u)=-Q\Sin_q(u);\quad
\cos_{q^{-1}}(u)=Q\Cos_q(u);\hfill\cr
\rlap{(2.10)}
\hfill \tan_{q^{-1}}(u)=-Q\tan_q(u);\qquad
\sec_{q^{-1}}(u)=Q\Sec_q(u).\hfill\cr}
$$
We end this Section with three technical lemmas that will be used in
Section~5 (resp. Section~6) to show how identity (1.12) (resp. (1.16)) in 
Theorem~1.2 (resp. in Theorem~1.3) can be obtained
from (1.11) (resp. from (1.15)). 

\proclaim Lemma 2.2. We have $D_{q^{-1}}Q=q\,QD_{q^{-1}}$,
so that for each $n\ge 1$
$$
(D_{q^{-1}}Q)^n=q^{(n+1)n/2}Q^nD_{q^{-1}}^n.\leqno(2.11)
$$

{\it Proof}.\quad
By (2.6) we have: $D_{q^{-1}}Q=U(I-Q^{-1})Q=UQ-U
=q\,QU-q\,QUQ=q\,QU(I-Q^{-1})=q\,QD_{q^{-1}}$. Identity (2.11) is then a simple consequence, as there are $(n+1)n/2$ transpositions
$QD_{q^{-1}}\leftrightarrow D_{q^{-1}}Q$ to be made to go from
$(QD_{q^{-1}})^n$ to $Q^nD_{q^{-1}}^n$.\qed

\proclaim Lemma 2.3. For each composition 
${\bf c}=(c_0,c_1,\ldots,c_m)$ of an integer $n\ge 1$,
we have
$$
Q^n\tan_q((q^{-1})^{\rho\,{\bf c}}u)=\tan_q(q^{\bf c}u).
\leqno(2.12)
$$

{\it Proof}.\quad
A simple verification:
$$
\displaylines{\quad
\tan_q((q^{-1})^{\rho\,{\bf c}}u)=\tan_q((q^{-1})^{c_m}u)
\tan_q((q^{-1})^{c_m+c_{m-1}}u)\hfill\cr
\kern6cm\times{}\cdots\times
\tan_q((q^{-1})^{c_m+\cdots+c_{1}}u).\hfill\cr}
$$
Hence, the left-hand side of (2.12) is equal to
$$
\leqalignno{\quad
Q^n\tan_q((q^{-1})&^{\rho\,{\bf c}}u)\cr
=&\tan_q(q^{n-c_m}u)
\tan_q(q^{n-c_m-c_{m-1}}u) \times \cdots\times
\tan_q(q^{n-c_m-\cdots-c_{1}}u)\cr
=&\tan_q(q^{c_0+\cdots+c_{m-1}}u)
\tan_q(q^{c_0+\cdots+c_{m-2}}u)
\times{}\cdots\times
\tan_q(q^{c_{0}}u)\hfill\cr
=&\tan_q(q^{\bf c}u).\qed\hfill\cr
}$$

\proclaim Lemma 2.4. We have:
$$
D_q^n\Sec_q(u)=(-1)^nq^{n(n-1)/2}Q^{n-1}D_{q^{-1}}^n\sec_{q^{-1}}(u).
$$

{\it Proof}.\quad
Just write
$$
\eqalignno{
D_q^n\Sec_q(u)&=(-D_{q^{-1}}Q)^n\,Q^{-1}\sec_{q^{-1}}(u)\qquad\qquad
&\hbox{[by (2.7) and (2.10)]}\cr
&=(-1)^n(D_{q^{-1}}Q)^{n-1}D_{q^{-1}}\sec_{q^{-1}}(u)\cr
&=(-1)^nq^{n(n-1)/2}Q^{n-1}D_{q^{-1}}^n\sec_{q^{-1}}(u).\qed
&\hbox{[by (2.11)]}\cr
}$$

\medskip
\centerline{\bf 3. Transformations on $t$-permutations} 

\medskip
Say that a $t$-permutation  
$w=(w_0,w_1,\ldots, w_m)$ is of the
{\it first} (resp. {\it second\/}) {\it kind}, if~$1$ appears (resp. does not
appear) as a one-letter factor among the~$w_i$'s. Each
set ${\cal T}_{n,m}$  can be
partitioned into two subsets
${\cal T}_{n,m}^{*}$ and ${}^*{\cal T}_{n,m}$, the former
one consisting of all permutations from ${\cal T}_{n,m}$
 of the first kind, the latter
one of those of the second kind. Let
$[m]:=\{1,2,\ldots,m\}$ and  for each integer~$y$ let 
$y^+:=y+1$
and $v^+:=y^+_1y^+_2\ldots y^+_m$ for each word
$v=y_1y_2\ldots y_m$, whose letters are integers.

For each pair $(m,n)$ such that
$0\leq m \leq n+1$ and $m\equiv n+1 \bmod 2$,
we construct two bijections 
$$
\Delta\kern-2pt^* :[m]\times {\cal
T}_{n,m}\rightarrow {\cal T}_{n+1,m+1}^*;\qquad
{}^*\kern-3pt\Delta :[m]\times {\cal T}_{n,m}\rightarrow 
{}^*{\cal T}_{n+1,m-1};\leqno(3.1)
$$
Let $w=(w_0,w_1,w_2,\ldots, w_m)$ belong to 
${\cal T}_{n,m}$ and $1\le i\le m$. The sequence
$$(w_0^+,\ldots, w_{i-2}^+,w_{i-1}^+,\,1\,, {w_i^+},
{w_{i+1}^+},\ldots, {w_m^+})\leqno(3.2)$$
is obviously  from ${\cal T}_{n+1,m+1}^{*}$; denote it by
$\Delta\kern-2pt{}^* (i,w)$. Next, the sequence
$$
(w_0^+,\ldots ,w_{i-2}^+,(w_{i-1}^+\,{1}\,
{w_i^+}),{w_{i+1}^+},\ldots, w_m^+)\leqno(3.3)
$$
is then from $^*{\cal T}_{n+1,m-1}$, and will be
denoted  by
${}^*\kern-3pt\Delta(i,w)$; its $i$-th factor
appears as
$w_{i-1}^+\, 1\; {w_i^+}$, as indicated in the inner 
parentheses in~(3.3). 
By the above two bijections, recurrence $(1.3)$ holds with the
interpretation $a(n,m)=\#{\cal T}_{n,m}$.

The next step is to study the actions of the transformations 
${}^*\kern-3pt\Delta$ and $\Delta\kern-2pt^*$
on the statistics ``ides,", ``imaj," ``min," ``inv" introduced in Subsection~1.4. Taking again the example used in that Subsection, namely the $t$-permutation
$w=({4}\,{5},\,11\,{1}\,{\bf3},\, {\bf 1\!0}\, {7}\,
{\bf9},\,{\bf6},\, {\bf8}\, {\bf2})$, whose underlying statistics are
$\ides w=6$, $\imaj w=38$,
$\inv w=27$ and
$\min w=1$, we get:

\goodbreak
$$\vbox{\halign{\hfil$#$&$\;=\;#$&$#$\hfil\quad&\ \hfil$#$\
\hfil\ &\ \hfil $#$\hfil\ &\ \hfil$#$\hfil\quad
&\ \hfil$#$\hfil\cr 
&\omit&&\ides&\imaj&\min&\inv\cr
w&&({4}\,
{5},\,11\,{1}\,{\bf3},\, {\bf 1\!0}\, {7}\,
{\bf9},\,{\bf6},\, {\bf8}\, {\bf2})
&6&38&1&27\cr
\noalign{\medskip}
\Delta\kern-2pt^*  (1,w)&&({5}\,
{6},\,{ 1},\,{1\!2}\, {2}\, {\bf4},\,
{\bf1\!1}\,{8}\,{\bf1\!0},\, {\bf7},\,{\bf9}\,
{\bf3})&6&44&1&29\cr 
\Delta\kern-2pt^*  (2,w)&&({5}\,
{6},\,{1\!2}\, {2}\,  {\bf4},\, {\bf1},\,
{\bf11}\,{8}\,{\bf1\!0},\, {\bf7},\,{\bf9}\,
{\bf3})&7&45&2&32\cr
\Delta\kern-2pt^*  (3,w)&&({5}\,
{6},\,{1\!2}\, {2}\, {\bf4},\,
{\bf1\!1}\,{8}\,{\bf1\!0},\,
{\bf1},\, {\bf7},\,{\bf9}\, {\bf3})&7&45&3&35\cr
\Delta\kern-2pt^*  (4,w)&&({5}\,
{6},\,{1\!2}\, {2}\, {\bf4},\,
{\bf1\!1}\,{8}\,{\bf1\!0},\,
{\bf7},\, {\bf1},\,{\bf9}\, {\bf3})&7&45&4&36\cr
\noalign{\medskip}
{}^*\kern-3pt\Delta(1,w)&&({ 5}\,
{6}\, {1}\,{1\!2}\, {2}\,
{\bf4},\, {\bf1\!1}\,{8}\,{\bf1\!0},\, {\bf7},\,{\bf9}\,
{\bf3})&6&44&0&29\cr 
{}^*\kern-3pt\Delta(2,w)&&({5}\,
{6},\,{1\!2}\, {2}\, {\bf4}\ {\bf1}\,
{\bf1\!1}\,{8}\,{\bf1\!0},\, {\bf7},\,{\bf9}\,
{\bf3})&7&45&1&32\cr
{}^*\kern-3pt\Delta (3,w)&&({5}\,
{6},\,{1\!2}\, {2}\, {\bf4},\,
{\bf1\!1}\,{8}\,{\bf1\!0}\, {\bf1}\,
{\bf7},\,{\bf9}\, {\bf3})&7&45&2&35\cr
{}^*\kern-3pt\Delta(4,w)&&({5}\,
{6},\,{1\!2}\, {2}\, {\bf4},\,
{\bf1\!1}\,{8}\,{\bf1\!0},\,
{\bf7}\,{\bf1}\,{\bf9}\, {\bf3})&7&45&3&36\cr
}}
$$

The proof of the next theorem is a simple verification and
will not be reproduced here.

\proclaim Theorem 3.1.  Let $w$ 
be a $t$-permutation and let $\Iligne
w:=\{j_1<j_2<\cdots<j_r\}$. 
Furthermore, let $\Delta^*(i,w)$ and ${}^*\Delta(i,w)$ be defined
as in $(3.1)$ - $(3.3)$. 
Then,
$$\displaylines{\quad
\Iligne \Delta\kern-2pt^* (i,w)\hfill\cr
\noalign{\vskip-6pt}
\hfill{}
=\Iligne {}^*\kern-3pt\Delta
(i,w)=\cases{\{(j_1+1), (j_2+1),\ldots,(j_r+1)\},&if
$i\le \min w$;\cr
\{1,(j_1+1), (j_2+1),\ldots,(j_r+1)\},&if
$\min w<i$.\cr}\quad\cr
\noalign{\hbox{Hence,}}
\noalign{\vskip -5pt}
\eqalign{\ides\Delta\kern-2pt^*
(i,w)&=\ides{}^*\kern-3pt\Delta (i,w)=
\cases{\ides w,&if $i\le \min w$;\cr
1+\ides w,&if $\min w<i$.\cr}\cr
\imaj\Delta\kern-2pt^*
(i,w)&=\imaj{}^*\kern-3pt\Delta (i,w)=
\cases{\ides w+\imaj w,&if $i\le \min w$;\cr
1+\ides w+\imaj w,&if $\min w<i$;\cr}\cr
}\cr
\noalign{\hbox{Furthermore,}}
\inv\Delta\kern-2pt^* (i,w)=\inv{}^*\kern-3pt\Delta (i,w) =\inv w+\lambda(w_0w_1\cdots w_{i-1}),\cr
\min(\Delta\kern-2pt^*{}(i,w))=i,\quad \min(^*\kern-3pt\Delta 
(i,w))=i-1,\ {\rm for}\ i\ge 1.
\cr }
$$

\medskip
\centerline{\bf 4. Proof of Theorem 1.1} 

\medskip
Two steps are needed to achieve the proof: first, the derivation of a {\it recurrence relation} for the polynomials $A_{n,k,a,b}(q)$, presented in the next theorem, then, an explicit {\it algorithm} for calculating them, described in Lemma~4.2.

\proclaim Theorem 4.1 {\rm (Recurrence for
the polynomials $A_{n,k,a,b}(q)$)}.
With $m:=a+b$ and $m':=a'+b'$ we have:
$$
\displaylines{(4.1)\ 
A_{n+1,k',a',b'}(q)\!=\!q^{k'}\Bigl(\,
\sum_{0\le a\le a'-1\le m'-2}\kern-25pt 
A_{n,k'-1,a,m'-1-a}(q)
+\kern-20pt \sum_{1\le a'\le a\le m'-1}\kern-20pt 
A_{n,k',a,m'-1-a}(q)\hfill\cr
\hfill{}+\sum_{0\le a\le a'}
A_{n,k'-1,a,m'+1-a}(q)
+\kern-15pt \sum_{a'+1\le a\le m'+1}
\kern-15pt A_{n,k',a,m'+1-a}(q)\Bigr),\cr
}
$$
valid for $n\ge 0$ with the initial condition:
$A_{0,k,a,b}(q)=\delta_{k,0}\delta_{a,1}\delta_{b,0}$.
\goodbreak

{\it Proof}.\quad 
Let
$w'=(w'_0,w'_1,\ldots,w'_{m'})\in{\cal T}_{n+1,k',a',b'}$.
When~$w'$ is of the first kind, the inverse
${\Delta\kern-2pt^*}^{-1}(w')$ is obtained by  deleting the
unique $w'_{a'}$ equal to~1 and subtracting~1 from all
the other letters of the components $w'_j$
$(j=0,\dots,a'-1,a'+1,\ldots,m')$. When~$w'$ is of the
second kind, the inverse
${}^*\kern-3pt\Delta^{-1}(w')$ is obtained by deleting~1
from the component $w'_{a'}=u'_{a'}1v'_{a'}$ and
subtracting~1 from all the letters of the components
$w'_0$, \dots~,
$w'_{{a'}-1}$,  $u'_{a'}$, $v'_{a'}$, $w'_{{a'}+1}$, \dots~,
$w'_{m'}$. The $t$-permutations $w'$ from
${\cal T}_{n+1,k',a',b'}$ fall into four categories:

\smallskip
(1) $w'$ is of the first kind and 2 is {\it to the left of}~1
in the product $w'_0w'_1\cdots w'_{m'}$; then
$1\le a'\le m'-1$, because, when $n+1\ge 2$, the
components $w'_0$ and $w'_{m'}$ are, either empty, or
their lengths are at least equal to~2. Hence,
${\Delta\kern-2pt^*}^{-1}(w')=(a,w)$ with
$w\in {\cal T}_{n,k'-1,a,m'-1-a}$ and $0\le a\le a'-1
\le m'-2$, as $\ides w=\ides w'-1$. Also, $\imaj
w'=1+\ides w+\imaj w=k'+\imaj w$. [By Theorem 3.1]

\smallskip
(2) $w'$ is of the first kind and 2 is {\it to the right
of}~1 in  $w'_0w'_1\cdots w'_{m'}$; for an analogous
reason as in case~(1) we have: 
${\Delta\kern-2pt^*}^{-1}(w')=(a,w)$ with
$w\in {\cal T}_{n,k',a,m'-1-a}$ and $a'\le a\le m'-1$. 
In this case $\imaj w'=\ides w+\imaj w=k'+\imaj w$.

\smallskip
(3) $w'$ is of the second kind and 2 is {\it to the left
of}~1 in $w'_0w'_1\cdots w'_{m'}$, so that~2 can
still belong to $w'_{a'}$, or to any one of the components
$w'_0$, \dots~, $w'_{a'-1}$. Hence,
${\Delta\kern-2pt^*}^{-1}(w')=(a,w)$
with $w\in {\cal T}_{n,k'-1,a,m'+1-a}$ and
$0\le  a\le a'$, for the number of components has
increased by~1. Again, $\imaj w'=k'+\imaj w$.

\smallskip
(4) $w'$ is of the second kind and 2 is {\it to the right
of}~1 in $w'_0w'_1\cdots w'_{m'}$, so that~2 can
still belong to $w'_{a'}$, or to any one of the components
$w'_{a'+1}$, \dots~, $w'_{m'}$. Hence,
${\Delta\kern-2pt^*}^{-1}(w')=(a,w)$
with $w\in {\cal T}_{n,k',a,m'+1-a}$ and
$a'+1\le  a\le m'+1'$. Also, $\imaj w'=k'+\imaj w$.

\smallskip
Thus, identity~(4.1) holds.\qed
\medskip

A consequence of the combinatorial interpretation is the 
symmetry property
$$
q^{n(n-1)/2}A_{n,n-1-k,b,a}(q^{-1})
=A_{n,k,a,b}(q),
$$
whose proof is easy and will be omitted.
\medskip
\proclaim Lemma 4.2.
Let $[k,a,b]:=(\tan_q(q^{k+1}u))^b\,(\tan_q(q^ku))^a$. Then,
$$\displaylines{
D_q[k,a,b]=q^k
\kern-10pt \sum_{0\le i\le a-1}\kern-8pt 
[k,i,a+b-1-i]+q^{k+1}\kern-10pt
\sum_{0\le i\le b-1}\kern-8pt
[k+1,a+i,b-1-i]\hfill\cr
\hfill {}+q^k
\kern-3pt \sum_{1\le i\le a}\kern-2pt 
[k,i,a+b+1-i]+q^{k+1}\kern-6pt
\sum_{1\le i\le b}\kern-2pt
[k+1,a+i,b+1-i].\cr
}
$$

\medskip
{\it Proof}.
For taking the $q$-derivative of a product $f_1(u)f_2(u)\cdots f_n(u)$ of 
$q$-series we use the formula:
$$\leqalignno{\qquad
D_q \prod_{1\le i\le n}f_i(u)
&=\sum_{1\le i\le n}
f_1(u)\cdots f_{i-1}(u)\,(D_qf_i(u))\,
f_{i+1}(qu)\cdots f_n(qu).&(4.2)\cr
}
$$

In particular, 
$$\displaylines{
\ D_q   (\tan_q(q^{k+1}u))^b\,(\tan_q(q^ku))^a\hfill \cr
\quad{}=\sum_{0\le i\le a-1}  (\tan_q(q^{k+1}u))^b\,
	(\tan_q(q^ku))^i (D_q \tan_q(q^ku)) (\tan_q(q^{k+1}u))^{a-1-i}\hfill\cr
\qquad{}+\sum_{0\le j\le b-1} \kern-8pt 
(\tan_q(q^{k+1}u))^j
(D_q\tan_q(q^{k+1}u))
(\tan_q(q^{k+2}u))^{b-1-j}
(\tan_q(q^{k+1}u))^a\hfill\cr
\quad{}=q^k\sum_{0\le i\le a-1}  (\tan_q(q^{k+1}u))^{a+b-1-i}\, (\tan_q(q^ku))^i\hfill \cr
\qquad{}+q^k\sum_{0\le i\le a-1}  (\tan_q(q^{k+1}u))^{a+b-i}\, (\tan_q(q^ku))^{i+1}\hfill \cr
\qquad{}+q^{k+1}\sum_{0\le j\le b-1}   (\tan_q(q^{k+2}u))^{b-1-j}(\tan_q(q^{k+1}u))^{a+j}  \hfill\cr
\qquad{}+q^{k+1}\sum_{0\le j\le b-1}   (\tan_q(q^{k+2}u))^{b-j}(\tan_q(q^{k+1}u))^{a+j+1} .\qed\hfill\cr
}
$$

With the notation $[k,a,b]$ we have $[0,1,0]=\tan_q(u)$, and identity (1.9)
can be rewritten 
$$
D_q^n\tan_q(u)=\sum_{k,a,b}A_{n,k,a,b}(q)\,[k,a,b],
$$
so that Lemma 4.2 can be used to calculate  the polynomials $A_{n,k,a,b}(q)$
by iteration :
$$\leqalignno{
D_q[0,1,0]&\ (=D_q\tan_q(u))=[0,0,0]+[0,1,1],\cr
\noalign{\hbox{\qquad so that $A_{1,0,0,0}(q)=A_{1,0,1,1}(q)=1$;}}
D_q^2[0,1,0]&\ (=D_q^2\tan_q(u))=D_q[0,1,1]=
[0,0,1]
+q[1,1,0]+[0,1,2]+q[1,2,1], \cr
\noalign{\hbox{\qquad so that $A_{2,0,0,1}(q)=1$, $A_{2,1,1,0}(q)=q$, 
$A_{2,0,1,2}(q)=1$, $A_{2,1,2,1}(q)=q$;}}
D_q^3[0,1,0]&\ (=D_q^3\tan_q(u))=
(q+q^2)[1,0,0]\cr
&\quad{}+[0,0,2]+q^2[1,0,2]+(2q+2q^2)[1,1,1]
+q[1,2,0]+q^3[2,2,0]\cr
&\quad{}+[0,1,3]+q^2[1,1,3]+(q+q^2)[1,2,2]
+q[1,3,1]+q^3[2,3,1],\cr
\noalign{\hbox{\qquad so that $A_{3,1,0,0}(q)=q+q^2$, $A_{3,0,0,2}(q)=1$, etc.}}
}
$$
The polynomials $A_{n,k,a,b}(q)$ in Table 2 (Section 11)
have been obtained using the previous calculation.

\medskip

{\it Proof of Theorem 1.1}.
By induction. Assume that (1.9) is true for $n$.
$$
\leqalignno{
D_q^{n+1}[0,1,0]
&=D_qD_q^{n}[0,1,0]\cr
&=D_q\sum_{k,a,b}A_{n,k,a,b}(q)[k,a,b]
=\sum_{k,a,b}A_{n,k,a,b}(q)D_q[k,a,b]\cr
&=\sum_{k,a,b}A_{n,k,a,b}(q)\Bigl(q^k
 \sum_{0\le i\le a-1}\kern-8pt 
[k,i,a+b-1-i]\cr
&\quad{}+q^{k+1}\kern-10pt
\sum_{0\le i\le b-1}\kern-8pt
[k+1,a+i,b-1-i]\cr
&\quad{}+q^k
\kern-3pt \sum_{1\le i\le a}\kern-2pt 
[k,i,a+b+1-i]+q^{k+1}\kern-6pt
\sum_{1\le i\le b}\kern-2pt
[k+1,a+i,b+1-i]\Bigr).\cr}
$$
We calculate the contribution
of each sum to the triple $[k',a',b']$. For the first sum we
have $[k,i,a+b-1-i]=[k',a',b']$ $(0\le i\le a-1)$ if and
only if $k=k'$, $i=a'$, $a+b-1=a'+b'$, $a'+1\le a\le
a'+b'+1$, so that the contribution is
$$
q^{k'}\sum_{a'+1\le a\le a'+b'+1}
A_{n,k',a,a'+b'+1-a}(q).
$$
For the second sum we have $[k+1,a+i,b-1-i]=[k',a',b']$
$(0\le i\le b-1)$ if and only if $k=k'-1$, $a+i=a'$, 
$b-1-i=b'$, $0\le a'-a\le b-1$, that is,
$k=k'-1$, $i=a'-a$, $a+b=a'+b'+1$, $0\le a\le a'$, so
that the contribution is
$$
q^{k'}\sum_{0\le a\le a'}
A_{n,k'-1,a,a'+b'+1-a}(q).
$$
with the convention that  $A_{n, -1, a, a'+b'+1-a}(q)=0$.
For the third sum we have $[k,i,a+b+1-i]=[k',a',b']$
$(1\le i\le a)$ if and only if
$k=k'$, $i=a'$, $a+b+1=a'+b'$, $1\le a'\le a\le
a+b=a'+b'-1$, so that the contribution is
$$
q^{k'}\sum_{1\le a'\le a\le a'+b'-1}
A_{n,k',a,a'+b'-1-a}(q).
$$
For the fourth sum we have $[k+1,a+i,b+1-i]=[k',a',b']$ 
$(1\le i\le b)$ if and only if $k=k'-1$, $a+i=a'$,
$b+1-i=b'$, $1\le i\le b$, that is, $k=k'-1$, $i=a'-a$,
$a+b+1=a'+b'$, $0\le a\le a'-1=a+i-1\le a+b-1\le
a'+b'-2$.

The contribution is then
$$q^{k'}
\sum_{0\le a\le a'-1\le a'+b'-2}
A_{n,k'-1,a,a'+b'-1-a}(q).
$$
By Theorem 4.1 we have
$$D_q^{n+1}[0,1,0]
=\sum_{k',a',b'}A_{n+1,k',a',b'}[k',a',b']. \qed$$
\medskip
\centerline{\bf 5. Proof of Theorem 1.2} 

\medskip
As for Theorem 1.1, two steps are used to complete the proof of
Theorem~1.2: first, the derivation of a recurrence relation, then, the construction of an explicit algorithm.

\goodbreak
\proclaim Theorem 5.1 {\rm (Recurrence relation for the
polynomials $B_{n,k,a,b}(q)$)}. With $m:=a+b$ and
$m':=a'+b'$ we have:
$$
\displaylines{(5.1)\ 
B_{n+1,k',a',b'}(q)=q^{k'}\Bigl(\;
\sum_{0\le a\le a'-1}\kern-10pt 
B_{n,k'-1,a,m'-1-a}(q)
+\kern-20pt \sum_{1\le a'\le a\le m'-1}\kern-20pt 
B_{n,k',a,m'-1-a}(q)\hfill\cr
\hfill{}+\sum_{0\le a\le a'}
B_{n,k'-1,a,m'+1-a}(q)
+\kern-15pt \sum_{a'+1\le a\le m'+1}
\kern-15pt B_{n,k',a,m'+1-a}(q)\Bigr),\cr
}
$$
valid for $n\ge 0$ with the initial condition:
$B_{0,k,a,b}(q)=\delta_{k,-1}\delta_{a,0}\delta_{b,0}$.

{\it Proof}.\quad 
Let
$w'=(w'_0,w'_1,\ldots,w'_{m'}, \epsilon)\in{\cal T}'_{n+1,k',a',b'+1}$.
What
has been said for recurrence~(4.1) can be reproduced, except
for the first sum, when $w'$ is of the first kind and 2 is
{\it to the left of}~1 in the product $w'_0w'_1\cdots
w'_{m'}$; now, $1\le a'\le m'$ and not $m'-1$. Hence,
${\Delta\kern-2pt^*}^{-1}(w')=(a,w)$ with
$w\in {\cal T'}_{n,k'-1,a,m'-a}$ and $0\le a\le a'-1
\le m'-1$.
Thus, identity~(5.1) holds.\qed

\proclaim Lemma 5.2.
Let $\langle k,a,b\rangle:=(\tan_q(q^{k+1}u))^b\,\sec_q(q^{k+1}u) (\tan_q(q^ku))^a$. Then
$$\displaylines{
\rlap{\hbox{\rm (5.2)}}\hfill 
D_q\langle k,a,b\rangle :=q^k\sum_{0\le i\le
a-1}\langle k,i,a+b-1-i\rangle +q^{k}\sum_{1\le i\le
a}\langle k,i,a+b+1-i\rangle\cr
\qquad\quad{}+q^{k+1}\sum_{0\le i\le b-1} \langle k+1,a+i,b-1-i\rangle +q^{k+1}\sum_{1\le i\le
b+1}\kern-10pt 
\langle k+1,a+i,b+1-i\rangle.\cr
}
$$

{\it Proof}.
By using (4.2) we derive 
$$\displaylines{
\ D_q   (\tan_q(q^{k+1}u))^b\,\sec_q(q^{k+1}u)\,(\tan_q(q^ku))^a\hfill \cr
\quad{}=\sum_{0\le i\le a-1}  (\tan_q(q^{k+1}u))^{a+b-1-i}\,\sec_q(q^{k+1}u) 
	(\tan_q(q^ku))^i (D_q \tan_q(q^ku)) \hfill\cr
\qquad{}+\sum_{0\le j\le b-1} \kern-8pt 
(\tan_q(q^{k+1}u))^{a+j}
(D_q\tan_q(q^{k+1}u))
(\tan_q(q^{k+2}u))^{b-1-j}
\sec_q(q^{k+2}u)\hfill\cr
\qquad{}+  (\tan_q(q^{k+1}u))^b\,D_q(\sec_q(q^{k+1}u))\,(\tan_q(q^ku))^a\hfill 
\cr
\quad{}=q^k\sum_{0\le i\le a-1}  (\tan_q(q^{k+1}u))^{a+b-1-i}\, \sec_q(q^{k+1}u)(\tan_q(q^ku))^i\hfill \cr
\qquad{}+q^k\sum_{0\le i\le a-1}  (\tan_q(q^{k+1}u))^{a+b-i}\, \sec_q(q^{k+1}u)(\tan_q(q^ku))^{i+1}\hfill \cr
\qquad{}+q^{k+1}\sum_{0\le j\le b-1}   (\tan_q(q^{k+2}u))^{b-1-j}\sec_q(q^{k+2}u) (\tan_q(q^{k+1}u))^{a+j} \hfill\cr
\qquad{}+q^{k+1}\sum_{0\le j\le b-1}   (\tan_q(q^{k+2}u))^{b-j}\sec_q(q^{k+2}u)(\tan_q(q^{k+1}u))^{a+j+1} \hfill\cr
\qquad{}+  q^{k+1}(\tan_q(q^{k+1}u))^{a+b+1}\,\sec_q(q^{k+2}u).\qed\hfill 
\cr
}
$$

Lemma 5.2 provides a way to calculate the 
polynomials $B_{n, k,a,b}(q)$. As $\langle -1, 0,0\rangle=\sec_q(u)$,
we get
\medskip
\noindent
$D_q\langle -1,0,0\rangle=\langle 0,1,0\rangle$, 
so that $B_{1,0,1,0}(q)=1$;
\smallskip
\noindent
$D_q^2\langle -1,0,0\rangle=D_q\langle 0,1,0\rangle=
\langle 0,0,0\rangle
+\langle 0,1,1\rangle+q\langle 1,2,0\rangle$, 
\smallskip
so that $B_{2,0,0,0}(q)=B_{2,0,1,1}(q)=1$, 
$B_{2,1,2,0}(q)=q$;
\smallskip
\noindent
$\eqalign{D_q^3\langle -1,0,0\rangle&=
D_q\bigl(\langle 0,0,0\rangle
+\langle 0,1,1\rangle+q\langle 1,2,0\rangle\bigr)\cr
&=q\langle 1,1,0\rangle\cr
&\ {}+\bigl(\langle 0,0,1\rangle+\langle 0,1,2\rangle
+q\langle 1,1,0\rangle+q\langle 1,2,1\rangle+q\langle 1,3,0\rangle\bigr)\cr
&\ {}+\bigl(q^2\langle 1,0,1\rangle\!+\!q^2\langle 1,1,0\rangle
\!+\!q^2\langle 1,1,2\rangle\!+\!q^2\langle 1,2,1\rangle
\!+\!q^3\langle 2,3,0\rangle\bigr)\cr
&=\langle 0,0,1\rangle+q^2\langle 1,0,1\rangle+(2q+q^2)\langle 1,1,0\rangle\cr
&\ {}+\!\langle 0,1,2\rangle\!+\!q^2\langle 1,1,2\rangle
\!+\!(q+q^2)\langle 1,2,1\rangle\!+\!q\langle 1,3,0\rangle\!+\!q^3\langle 2,3,0\rangle,\cr
}
$
\smallskip
so that $B_{3,0,0,1}(q)=1$, $B_{3,1,0,1}(q)=q^2$, 
$B_{3,1,1,0}(q)=2q+q^2$, $B_{3,0,1,2}(q)=1$,
$B_{3,1,1,2}(q)=q^2$, $B_{3,1,2,1}(q)=q+q^2$,
$B_{3,1,3,0}(q)=q$, $B_{3,2,3,0}(q)=q^3$.
\medskip
The polynomials $B_{n,k,a,b}(q)$ in Table 2 have been calculated by means of
the previous algorithm.
\medskip

{\it Proof of (1.11)}.
By induction. Assume that (1.11) is true for $n$. Then,
$$
\displaylines{
D_q^{n+1}\langle -1,0,0\rangle 
= D_q D_q^{n}\langle -1,0,0\rangle 
\hfill\cr
\kern1cm{}=D_q\sum_{k,a,b}B_{n,k,a,b}(q)\langle k,a,b\rangle 
=\sum_{k,a,b}B_{n,k,a,b}(q)D_q\langle k,a,b\rangle \hfill\cr
\kern1cm{}=\sum_{k,a,b}B_{n,k,a,b}(q)\Bigl(q^k
 \sum_{0\le i\le a-1}\kern-8pt 
\langle k,i,a+b-1-i\rangle \hfill\cr
\kern1cm{}\quad{}+q^{k+1}\kern-10pt
\sum_{0\le i\le b-1}\kern-8pt
\langle k+1,a+i,b-1-i\rangle \hfill\cr
\kern1cm{}\quad{}+q^k
\kern-3pt \sum_{1\le i\le a}\kern-2pt 
\langle k,i,a+b+1-i\rangle +q^{k+1}\kern-6pt
\sum_{1\le i\le b+1}\kern-2pt
\langle k+1,a+i,b+1-i\rangle \Bigr).\hfill\cr}
$$
 We calculate the contribution
of each sum to the triple $\langle k',a',b'\rangle $. For the first
three sums we can simply reproduce the arguments
developed in the proof of Theorem~1.1. Only the fourth
sum is to be checked. This time, the double
inequality $1\le i\le b+1$ prevails, instead of $1\le i\le
b$. This leads to the sequence $0\le a\le a'-1=a+i-1\le
a+(b+1)-1\le a+b\le a'+b'-1$. Hence, the contribution is 
$$\displaylines{q^{k'}
\sum_{0\le a\le a'-1\le a'+b'-1}
B_{n,k'-1,a,a'+b'-1-a}(q).\cr
\noalign{\hbox{By Theorem 5.1 we have}}
 D_q^{n+1}\langle -1,0,0\rangle 
= \sum_{k',a',b'}B_{n+1,k',a',b'}(q)\langle k',a',b'\rangle .
\qed\cr}
$$

{\it Proof of (1.12)}.
Rewrite identity (1.11) taking (2.10) into account:
$$\displaylines{
D_{q^{-1}}^n\sec_{q^{-1}}(u)\hfill\cr
\quad{}=
\sum_{k,a,b}B_{n,k,a,b}(q^{-1})(\tan_{q^{-1}}(q^{-k-1}u))^b\,\sec_{q^{-1}}(q^{-k-1}u)(\tan_{q^{-1}}(q^{-k}u))^a\hfill\cr
\quad{}=
\sum_{k,a,b}B_{n,k,a,b}(q^{-1})(-1)^b(\tan_{q}(q^{-k}u))^b\,\Sec_{q}(q^{-k}u)(-1)^a(\tan_{q}(q^{-k+1}u))^a.\hfill\cr
}
$$
Hence, as $0\le k\le n-1$, 
$$\displaylines{
Q^{n-1}D_{q^{-1}}^n\sec_{q^{-1}}(u)\hfill\cr
\quad{}=
\sum_{k,a,b}B_{n,k,a,b}(q^{-1})(-1)^{a+b}(\tan_{q}(q^{n-1-k}u))^b
\Sec_{q}(q^{n-1-k}u)(\tan_{q}(q^{n-k}u))^a\hfill\cr
\quad{}=
\sum_{k,a,b}B_{n,n\!-\!1\!-k,a,b}(q^{-1})(-1)^{a+b}(\tan_{q}(q^{k}u))^b
\Sec_{q}(q^{k}u)(\tan_{q}(q^{k+1}u))^a.\hfill\cr}
$$
Finally, as $n$ and $a+b$ are of the same parity, we have by Lemma 2.4
$$\displaylines{
D_q^n\Sec_q(u)=(-1)^nq^{n(n-1)/2}Q^{n-1}D_{q^{-1}}^n\sec_{q^{-1}}(u),\hfill\cr
\quad {}=
\sum_{k,a,b}q^{n(n-1)/2}B_{n,n\!-\!1\!-k,a,b}(q^{-1})(\tan_{q}(q^{k}u))^b
\Sec_{q}(q^{k}u)(\tan_{q}(q^{k+1}u))^a.\hfill\cr
}$$
This establishes identity (1.12). \qed

\bigskip
\centerline{\bf 6. Proof of Theorem 1.3}  
\medskip
Let ${\bf c}'=(c'_0,c'_1,\ldots,c'_{m'})\in\Theta_{n+1}$ and 
$0\leq i\leq m'$. If $c'_0\ge 2$, 
and if $1\le 2j+1\le c'_0$, let
$$
{\bf c}'[0,2j+1]:=(2j,c'_0-(2j+1),c'_{1},
\ldots,c'_{m'}). \leqno(6.1)
$$
If $1\le i\le m'$ and $2\le 2j\le c'_i$, let
$$
{\bf c}'[i,2j]:=(c'_0,\ldots,c'_{i-1},2j-1,c'_i-2j,c'_{i+1},
\ldots,c'_{m'}). \leqno(6.2)
$$
Finally, if $1\le i\le m'-1$ and $c'_i=1$, let
$$
{\bf c'}[i,1]:=
(c'_0,\ldots,c'_{i-1},c'_{i+1},\ldots,c'_{m'}). \leqno(6.3)
$$
\goodbreak

\proclaim Theorem 6.1 {\rm (Recurrence
relation for the polynomial $A_{n,{\bf c}}(q)$)}. With
${\bf c}'=(c'_0,c'_1,\ldots, c'_{m'})$ we have:
$$
\leqalignno{
A_{n+1,{\bf c}'}(q)&=\sum_{0\le 2j+1\le c'_0}
q^{2j}A_{n,{\bf c}'[0,2j+1]}(q)&(6.4)\cr
&\kern1cm{}+\sum_{1\le i\le m'}
q^{c'_0+\cdots+c'_{i-1}}
\sum_{2\le 2j\le c'_i}q^{2j-1}A_{n,{\bf c}'[i,2j]}(q) \cr
&\kern1cm{}+\sum_{1\le i\le
m'-1}\chi(c'_i=1)q^{c'_0+\cdots+c'_{i-1}}A_{n,{\bf
c}'[i,1]}(q),\cr}
$$
where $\chi(c'_i=1)$ is equal to 1 if $c'_i=1$ holds and $0$ otherwise.

{\it Proof}.\quad
Let $w'=(w'_0,w'_1,\ldots,w'_{m'})$ be a $t$-permutation
from ${\cal T}_{n+1}$ such that $\Lambda w'={\bf
c}'=(c'_0,c'_1,\ldots, c'_{m'})$. 
Two cases are to consider: (i) 1 belongs to the
component~$w_i$ $(0\le i\le m')$ having at least two
letters: write $w\in{\cal T}_{n+1,{\bf c}',2}$; (ii) the
component $w_i$ is equal to the
one-letter~1: write $w\in{\cal T}_{n+1,{\bf c}',1}$. For
each word~$v$, whose letters are integers, let $v^-$
designate the word obtained from~$v$ by subtracting~1
from each of its letters.  By convention,
$\epsilon^-:=\epsilon$.

Case (i): we may write: $w'_i=u'_i1v'_i$,
where at least one of the factors $u'_i$, $v'_i$ is
nonempty. 
To~$w'$ there corresponds a unique triple
$(w,i,\lambda(u'_i1))$ where
$w:=(w'{}^-_{\kern-3pt 0},\ldots,
w'{}^-_{\kern-3pt i-1},
u'{}^-_{\kern-3pt i},
v'{}^-_{\kern-3pt i},
w'{}^-_{\kern-3pt i+1},
w'{}^-_{\kern-3pt m'})$.
As $w'_0$ is rising alternating (when
nonempty) and the other components falling alternating,
the length $\lambda(u'_i1)$ of the word $u'_i1$ is odd,
say, $2j+1$, when $i=0$ and even, say, $2j$, when $i\ge
1$. 

Thus, the mapping $w\mapsto (w,i,\lambda(u'_i1))$ is a
bijection of ${\cal T}_{n+1,{\bf c}',2}$ to the set of
triples $(w,i,j)$, where  $0\le i\le m'$; $w\in {\cal
T}_n$, 
$\Lambda w$ is equal to ${\bf c}'[0,2j+1]$ when $i=0$
and to ${\bf c}'[i,2j]$ when $1\le i\le m'$;
 $0\le 2j+1\le c'_0$ when $i=0$ and $2\le 2j\le c'_i$
when $1\le i\le m'$. Moreover,
$$\inv w'=\cases{2j+\inv w,&when $i=0$;\cr
c'_0+\cdots+c'_{i-1}+2j-1+\inv w,&when $1\le i\le
m'$.\cr}
$$

Case (ii): the mapping $w'\mapsto (w,i)$ is a bijection of
${\cal T}_{n+1,{\bf c}',1}$ to the set of the pairs $(w,i)$,
where $w:=(w'{}^-_{\kern-3pt 0},\ldots,
w'{}^-_{\kern-3pt i-1},
w'{}^-_{\kern-3pt i+1},
w'{}^-_{\kern-3pt m'})$ and $1\le i\le m'-1$. Thus,
$\Lambda w={\bf c}'[i,1]$. Moreover,
$\inv w'=c'_0+\cdots+c'_{i-1}+\inv w$.

Accordingly, (6.4) holds.\qed
\medskip

Let ${\bf
c}=(c_0,c_1,\ldots,c_m)$ be a composition of a nonnegative
integer~$n$. If $1\le i\le m$,~let
$$
\leqalignno{
{}^{(i)}{\bf c}&:=(c_0,\ldots,c_{i-2},
c_{i-1}+1+c_i,c_{i+1},\ldots ,c_m);&(6.5)\cr
{\bf c}^{(i)}&:=(c_0,\ldots,
c_{i-1},1,c_{i},\cdots,c_m).&(6.6)\cr
}
$$
With the notations of (6.1)--(6.3) and (6.5)--(6.6) we then have:
$$
\leqalignno{
{}^{(1)}({\bf c}'[0,2j+1])&={\bf c}'.&(6.7)\cr
{}^{(i+1)}({\bf c}'[i,2j])&={\bf c}'.&(6.8)\cr
({\bf c'}[i,1])^{(i)}&={\bf c}'.&(6.9)\cr
}
$$

\proclaim Lemma 6.2.
Let ${\bf c}=(c_0,c_1,\ldots, c_m)\in\Theta_n$. Then, 
$$D_q\tan_q(q^{\bf c}u)
=\sum_{1\le i\le m} q^{c_0+\cdots+c_{i-1}}
\bigl(\tan_q(q^{{}^{(i)}{\bf c}}u)+\tan_q(q^{{\bf c}^{(i)}}u)\bigr).$$

{\it Proof}.
$$
\leqalignno{
D_q\tan_q(q^{\bf c}u)
&=D_q{\textstyle\prod\limits_{1\le i\le
m}}\tan_q(q^{c_0+\cdots+c_{i-1}}u)\cr
&=\sum_{1\le i\le m}\tan_q(q^{c_0}u)\cdots
\tan_q(q^{c_0+\cdots+c_{i-2}}u)\cr
&\qquad{}\times
q^{c_0+\cdots+c_{i-1}}
\bigl(1+\tan_q(q^{c_0+\cdots+c_{i-1}}u)
\tan_q(q^{c_0+\cdots +c_{i-1}+1}u)\bigr)\cr
&\qquad{}\times
\tan_q(q^{c_0+\cdots+c_{i-1}+1+c_i}u)
\cdots \tan_q(q^{c_0+\cdots+c_{i-1}+1+c_i
+\cdots+c_{m-1}}u)\cr
&=\sum_{1\le i\le m} q^{c_0+\cdots+c_{i-1}}
\bigl(\tan_q(q^{{}^{(i)}{\bf c}}u)+\tan_q(q^{{\bf c}^{(i)}}u)\bigr).\qed\cr
}
$$

In the next example, let ${\bf c}:=\tan_q(q^{\bf c}u)$ by convention, so that
$$D_q{\bf c}
=\sum_{1\le i\le m} q^{c_0+\cdots+c_{i-1}}
\bigl({{}^{(i)}{\bf c}}+{{\bf c}^{(i)}}\bigr).$$
We get
$$
\eqalignno{
D_q ((0,0))&:=(1)+(0,1,0)&\cr
D_q ^2((0,0))&:=D_q \bigl((1)+(0,1,0)\bigr)\cr
&=0+\bigl((2,0)+q(0,2)+(0,1,1,0)+q(0,1,1,0)\bigr)\cr
&=(2,0)+q(0,2)+(1+q)(0,1,1,0);\cr
D_q ^3((0,0))&:=D_q \bigl(
(2,0)+q(0,2)+(1+q)(0,1,1,0)\bigr)\cr
&=\bigl(q^2(3)+q^2(2,1,0)\bigr)+q\bigl((3)+(0,1,2)\bigr)\cr
&\qquad{}+(1+q)((2,1,0)+q(0,3,0)+q^2(0,1,2)\cr
&\qquad\qquad\qquad\quad {}+(0,1,1,1,0)
+q(0,1,1,1,1,0)+q^2(0,1,1,1,0)\bigr)\cr
&=(q+q^2)(3)+(1+q+q^2)(2,1,0)+(q+q^2+q^3)(0,1,2)\cr
&\qquad\qquad\qquad\quad{}
+(q+q^2)(0,3,0)+(1+q)(1+q+q^2)(0,1,1,1,0).\cr
}
$$
We can then calculate the polynomials (see also Table 3):
$A_{1,(1)}(q)=A_{1,(0,1,0)}(q)=1$;\hfil\break
$A_{2,(2,0)}(q)=1$, $A_{2,(0,2)}(q)=q$, 
$A_{2,(0,1,1,0)}(q)=1$, $A_{2,(0,1,1,0)}(q)=q$;\hfil\break
$A_{3,(3)}(q)=q+q^2$, 
$A_{3,(2,1,0)}(q)=1+q+q^2$, 
$A_{3,(0,1,2)}(q)=q+q^2+q^3$,
$A_{3,(0,3,0)}(q)=q+q^2$,
$A_{3,(0,1,1,1,0)}(q)=(1+q)(1+q+q^2)$.

\medskip
\goodbreak
{\it Proof of (1.14)}.\quad
Let ${\bf c}'=(c'_0,c'_1,\ldots, c'_{m'})$.
By (6.1) and (6.7) we have ${}^{(1)}{\bf c}={\bf c}'$ when
$c'_0\ge 2$, if and only if
${\bf c}={\bf c}'[0,2j+1]$ for some $j$ such that
$1\le 2j+1\le c'_0$. By (6.2) and (6.8) the relation
${}^{(i)}{\bf c}={\bf c}'$ holds for~$i$ such that $1\le
i\le m'$ and $c'_i\ge 2$, if and only if
${\bf c}={\bf c}'[i,2j]$ for some $j$ such that
$2\le 2j\le c'_i$. Finally, ${\bf c}^{(i)}={\bf c}'$ holds
for~$i$ such that $1\le i\le m'-1$ and $c'_i=1$, if and
only if
${\bf c}={\bf c}'[1,i]$.

Finally, (1.14) is proved by induction. Assume it is rue for $n$. Then,
$$\leqalignno{
D^{n+1}_q\tan_q(u)&=D_q D_q^n \tan_q(u)
=D_q \Bigl(\sum\limits_{{\bf c}\in \Theta_n} A_{n,{\bf
c}}(q)\,\tan_q(q^{\bf c}u)\Bigr)\cr 
&=\sum\limits_{{\bf c}\in \Theta_n} A_{n,{\bf
c}}(q)\,D_q \tan_q(q^{\bf c}u)\cr 
&=\sum\limits_{{\bf c}\in \Theta_n} A_{n,{\bf
c}}(q)\,
\sum_{1\le i\le m} q^{c_0+\cdots+c_{i-1}}
\bigl(\tan_q(q^{{}^{(i)}{\bf c}}u)+\tan_q(q^{{\bf c}^{(i)}}u)\bigr).
\cr 
&=\sum\limits_{{\bf c}'\in \Theta_{n+1}} 
\Bigl(\sum_{0\le 2j+1\le c'_0}
q^{2j}A_{n,{\bf c}'[0,2j+1]}(q)\cr
&\kern1cm{}+\sum_{1\le i\le m'}
q^{c'_0+\cdots+c'_{i-1}}
\sum_{2\le 2j\le c'_i}q^{2j-1}A_{n,{\bf c}'[i,2j]}(q) \cr
&\kern1cm{}+\sum_{1\le i\le
m'-1}\chi(c'_i=1)q^{c'_0+\cdots+c'_{i-1}}A_{n,{\bf
c}'[i,1]}(q)\Bigr) \, \tan_q(q^{{\bf c}'}u)\cr
&=\sum\limits_{{\bf c}'\in \Theta_{n+1}} 
\Bigl(A_{n+1, {\bf c}'}(q)  \Bigr) \, \tan_q(q^{{\bf c}'}u).
	\quad \hbox{\rm [By Theorem 6.1]}\qed\cr
}
$$
\medskip

For each ${\bf c}=(c_0,c_1,\ldots, c_m,0)\in\Theta_{n}^-$, let
$\Phi({\bf c}) := \tan_q(q^{{\bf c}^-}u) \sec_q(q^nu)$.

\proclaim Lemma 6.3.
Let ${\bf c}=(c_0,c_1,\ldots, c_m,0)\in\Theta_{n}^-$. Then, 
$$
D_q\Phi({\bf c})
=\sum_{1\le i\le m}\Bigl(q^{c_0+\cdots+c_{i-1}}
(\Phi({}^{(i)}{\bf c})+\Phi({\bf c}^{(i)})\Bigr) + q^n \Phi {\bf c}^{(m+1)}.
$$

{\it Proof}.
$$
\leqalignno{D_q\Phi({\bf c})
&=D_q\tan_q(q^{{\bf c}^-}u) \sec_q(q^nu)\cr
&=\sum_{1\le i\le m}\tan_q(q^{c_0}u)\cdots
\tan_q(q^{c_0+\cdots+c_{i-2}}u)\cr
&\qquad{}\times
q^{c_0+\cdots+c_{i-1}}
\bigl(1+\tan_q(q^{c_0+\cdots+c_{i-1}}u)
\tan_q(q^{c_0+\cdots +c_{i-1}+1}u)\bigr)\cr
&\qquad{}\times
\tan_q(q^{c_0+\cdots+c_{i-1}+1+c_i}u)
\cdots \tan_q(q^{c_0+\cdots+c_{i-1}+1+c_i
+\cdots+c_{m-1}}u)\cr
&\qquad{}\times \sec_q(q^{n+1}u)\cr
&\quad{}+q^n \tan_q(q^{{\bf c}^-}u) \tan_q(q^nu)\sec_q(q^{n+1}u)\cr
&=\sum_{1\le i\le m}\Bigl(q^{c_0+\cdots+c_{i-1}}
(\Phi({}^{(i)}{\bf c})+\Phi({\bf c}^{(i)})\Bigr) + q^n \Phi {\bf c}^{(m+1)}.\qed\cr
}
$$
\medskip

{\it Proof of (1.15)}.
Again (1.15) is proved by induction.  
Assume that (1.15) holds  for $n$. Then,
$$\leqalignno{
D^{n+1}_q\sec_q(u)&=D_q D_q^n \sec_q(u)
=D_q \Bigl(\sum\limits_{{\bf c}\in \Theta_n^-} A_{n,{\bf
c}}(q)\,\Phi({\bf c})\Bigr)\cr 
&=\sum\limits_{{\bf c}\in \Theta_n^-} A_{n,{\bf
c}}(q)\,D_q \Phi({\bf c})\cr 
&=\sum\limits_{{\bf c}\in \Theta_n^-} A_{n,{\bf
c}}(q)\,
\sum_{1\le i\le m}\Bigl(q^{c_0+\cdots+c_{i-1}}
(\Phi({}^{(i)}{\bf c})+\Phi({\bf c}^{(i)})\Bigr) + q^n \Phi {\bf c}^{(m+1)}
\cr
&=\sum\limits_{{\bf c}'\in \Theta_{n+1}^-} 
\Bigl(\sum_{0\le 2j+1\le c'_0}
q^{2j}A_{n,{\bf c}'[0,2j+1]}(q)\cr
&\kern1cm{}+\sum_{1\le i\le m'}
q^{c'_0+\cdots+c'_{i-1}}
\sum_{2\le 2j\le c'_i}q^{2j-1}A_{n,{\bf c}'[i,2j]}(q) \cr
&\kern1cm{}+\sum_{1\le i\le
m'-1}\chi(c'_i=1)q^{c'_0+\cdots+c'_{i-1}}A_{n,{\bf
c}'[i,1]}(q)\Bigr) \, \Phi({{\bf c}'})\cr
&=\sum\limits_{{\bf c}'\in \Theta_{n+1}^-} 
A_{n+1, {\bf c}'}(q)   \, \Phi({{\bf c}'}).
	\quad \hbox{\rm [By Theorem 6.1]}\qed\cr
}
$$
\medskip

{\it Proof of (1.16)}.
In our proof of (1.16) given next we again make a full use of the duality
derived in Section~2 between the $q$-series $\Sec_q(u)$, $\tan_q(u)$ 
and their analogs $\sec_{q^{-1}}(u)$, $\tan_{q^{-1}}(u)$. 
Since  $\tan_{q^{-1}}((q^{-1})^{{\bf c}^-}u)
=(-1)^{\mu {\bf c}}\tan_{q}((q^{-1})^{{\bf c}^-}qu) $, 
identity (1.15) can be rewritten 
$$\eqalignno{
D_{q^{-1}}^n\sec_{q^{-1}}(u)
&=\sum_{{{\bf c}\in \Theta_n^{-}}}
A_{n,{\bf c}}(q^{-1})\,\tan_{q^{-1}}((q^{-1})^{{\bf c}^-}u)\,\sec_{q^{-1}}((q^{-1})^n u)\cr
& =(-1)^{\mu {\bf c}}\sum_{{{\bf c}\in \Theta_n^{-}}}
A_{n,{\bf c}}(q^{-1})\,\tan_{q}((q^{-1})^{{\bf c}^-}qu)\,\Sec_{q}((q^{-1})^n qu).\cr
}
$$
By Lemma 2.3,
$$
Q^{n-1}D_{q^{-1}}^n\sec_{q^{-1}}(u)
	=(-1)^{\mu {\bf c}}
	\sum_{{{\bf c}\in \Theta_n^{-}}}\kern-7pt
A_{n,{\bf c}}(q^{-1})\tan_{q}(q^{\rho{\bf c}^-}u)\Sec_{q}(u).
$$
\goodbreak\noindent
By Lemma 2.4, and since $n$ and $\mu {\bf c}$ are of the same parity, we get
$$
\leqalignno{
D_q^n\Sec_q(u)
&=(-1)^nq^{n(n-1)/2}Q^{n-1}D_{q^{-1}}^n\sec_{q^{-1}}(u)\cr
&=\sum_{{\bf c}\in \Theta_n^{-}}
q^{n(n-1)/2}
A_{n,{\bf c}}(q^{-1})
\tan_{q}(q^{\rho{\bf c}^-}u)
\Sec_{q}(u).\hfill\cr
}
$$
This proves identity (1.16). \qed

\medskip
\goodbreak
%
\centerline{\bf 7. More on $q$-trigonometric functions} 

\medskip
Some parts of this section are of
semi-expository nature, although, to our knowledge, the
combinatorial properties of ``$\Sec_q(u)$" have not been
explicitly written down. Let
$$\leqalignno{
\tan_q(u)&=\sum_{n\ge 0}A_{2n+1}(q){u^{2n+1}\over
(q;q)_{2n+1}};&(7.1)\cr
\sec_q(u)&=\sum_{n\ge 0}A_{2n}(q){u^{2n}\over
(q;q)_{2n}};&(7.2)\cr
\Sec_q(u)&=\sum_{n\ge 0}A^{\Sec}_{2n}(q){u^{2n}\over
(q;q)_{2n}};&(7.3)\cr
} 
$$
be the $q$-expansions of the three series $\tan_q(u)$,
$\sec_q(u)$, $\Sec_q(u)$, respectively. 

Let ${N\brack M}_q:=(q;q)_N/((q;q)_M\,(q;q)_{N-M})$ 
$(0\le M\le N)$ be 
the Gaussian polynomial. 
Identities (2.3)-(2.5) yield 
$$
\leqalignno{
\quad A_{2n+1}(q)&=\sum_{0\le k\le n-1} {2n\brack 2k+1}
q^{2k+1}
A_{2k+1}(q)\,A_{2n-2k-1}(q),&(7.4)\cr
A_{2n}(q)&=\sum_{0\le k\le n-1}
{2n-1\brack 2k}_qq^{2k}
A_{2k}(q)\,A_{2n-2k-1}(q),&(7.5)\cr
A^{\Sec}_{2n}(q)&=\sum_{0\le k\le n-1}
{2n-1\brack 2k}_q A^{\Sec}_{2k}(q)\,q^{2n-2k-1}
A_{2n-2k-1}(q),&(7.6)\cr
}	
$$
for $n\geq 1$ with the initial conditions $A_1(q)=1$,
$A_0(q)=1$ and $A^{\Sec}_0(q)=1$.
Hence, the coefficients $A_n(q)$,
$A^{\Sec}_{2n}(q)$ $(n\ge 0)$ occurring in the $q$-expansions
of $\tan_q(u)$, $\sec_q(u)$, $\Sec_q(u)$ in
$(7.1)$--$(7.2)$ are polynomials with positive integral
coefficients. 
 
\medskip
The first values of the polynomials
$A_n(q)$  and $A^{\Sec}_{2n}(q)$ $(n\ge 0)$ can be calculated
by means of (7.4)---(7.6):

$A_1(q)\!=\!1$;\  $A_3(q)\!=\!q+q^2$;\ 
$A_5(q)\!=\!q^2+2q^3+3q^4+4q^5+3q^6+2q^7+q^8$;

\smallskip
$A_0(q)=A_2(q)=1$;\ 
$A_4(q)=q+2q^2+q^3+q^4$;

$A_6(q)=q^2+3q^3+5q^4+8q^5
+10q^6+10q^7+9q^8+7q^9+5q^{10}+2q^{11}+q^{12}$;

\smallskip
$A^{\Sec}_0(q)=A^{\Sec}_2(q)=q$;\ 
$A^{\Sec}_4(q)=q^2+q^3+2q^4+q^5$;

$A^{\Sec}_6(q)=q^3+2q^4+5q^5
+7q^6+9q^7+10q^8+10q^9+8q^{10}+5q^{11}+3q^{12}
+q^{13}$;

\medskip
We note that the two identities (7.5)
and (7.6) can be combined into the following single
formula, valid for
$n\ge 1$ with $A_0(q)=1$:
$$
A_n(q)=\sum_{0\le k\le \lfloor n/2\rfloor-1}
{n-1\brack 2k+1}_qq^{n-2k-2}A_{2k+1}(q)
A_{n-2k-2}(q).
$$
The polynomials $A_{2n+1}(q)$ (resp. $A_{2n}(q)$)
defined by (7.1) and (7.2) are usually called
the {\it $q$-tangent numbers} and {\it $q$-secant
numbers}, respectively. No traditional name exists for the
polynomials $A^{\Sec}_{2n}(q)$, as they are intimately related
to the $A_{2n+1}(q)$'s by identity (7.12).

\medskip
For each $n\ge 0$ the set of all rising (resp. falling)
alternating permutations $\sigma=\sigma(1)\sigma(2)
\cdots \sigma(n)$ of $1\,2\,\ldots\,n$ is denoted by
${\cal RA}_n$ (resp. ${\cal F\!A}_n$). As already
mentioned in Section~7, the old result by D\'esir\'e Andr\'e
[An81] asserts that both $\#{\cal RA}_{2n+1}$ and $\#{\cal
F\!A}_{2n+1}$  (resp. $\#{\cal RA}_{2n}$ and $\#{\cal
F\!A}_{2n}$) are equal to the tangent number
$T_{2n+1}$  (resp. secant number $E_{2n}$) occurring in
(1.1) and (1.2)
 
For the $q$-analog of this result we keep the
same combinatorial set-up, namely, ${\cal RA}_n$ and
${\cal F\!A}_n$ $(n\ge 0)$, but, as there are two
different $q$-secants, $\sec_q(u)$ and $\Sec_q(u)$, the
coefficients in their $q$-expansions will have different
combinatorial interpretations. For each
permutation $\sigma=\sigma(1)\sigma(2)\cdots\sigma(n)$
of $12\cdots n$ (not necesarily an alternating permutation),
let $\inv\sigma$ denote the traditional {\it number of
inversions} of~$\sigma$. 

\proclaim Theorem 7.1. Let $A_{2n+1}(q)$ (resp.
$A_{2n}(q)$, resp.
$A_{2n}^{\Sec}(q)$) be the coefficients in the $q$-expansion of  
$\tan_q(u)$ (resp. of $\sec_q(u)$, resp. of $\Sec_q(u)$), as
defined in
$(7.1)$--$(7.2)$. Then,
$$\leqalignno{
A_{2n+1}(q)&=\sum_{\sigma\in {\cal RA}_{2n+1}}
q^{\inv\sigma}=\sum_{\sigma\in {\cal F\!A}_{2n+1}}
q^{\inv\sigma};&(7.7)\cr
A_{2n}(q)&=\sum_{\sigma\in {\cal RA}_{2n}}
q^{\inv\sigma};&(7.8)\cr
A^{\Sec}_{2n}(q)&=\sum_{\sigma\in {\cal F\!A}_{2n}}
q^{\inv\sigma}.&(7.9)\cr}
$$

The proofs of (7.7) -- (7.9) are not reproduced here. It suffices to 
$q$-mimick Desir\'e Andr\'e's [An81] classical proof.

\medskip
The statistics ``Ligne" and ``imaj" have been defined
in Section~8. By means of the so-called ``second fundamental
transformation" (see, e.g., [Lo83, \S\kern2pt10.6], [Fo68],
[FS78]) we can construct a bijection $\Phi$ of the
group of all permutations onto itself with the property that
$$
\Ligne\sigma=\Ligne\Phi(\sigma)
\quad{\rm and}\quad
\inv\sigma=\imaj\Phi(\sigma).\leqno(7.10)
$$
Saying that a permutation~$\sigma$ is falling (resp.
rising) alternating is equivalent to saying that
$\Ligne\sigma=\{1,3,5,\ldots\,\}$ (resp.
$=\{2,4,6,\ldots\,\}$). Accordingly, we also have
$$\leqalignno{
A_{2n+1}(q)&=\sum_{\sigma\in {\cal RA}_{2n+1}}
q^{\imaj\sigma}=\sum_{\sigma\in {\cal F\!A}_{2n+1}}
q^{\imaj\sigma};\cr
A_{2n}(q)&=\sum_{\sigma\in {\cal RA}_{2n}}
q^{\imaj\sigma};\cr
A^{\Sec}_{2n}(q)&=\sum_{\sigma\in {\cal F\!A}_{2n}}
q^{\imaj\sigma}.\cr}
$$

For each permutation $\sigma=\sigma(1)\sigma(2)\cdots
\sigma(n)$ let ${\rho}$ (the {\it mirror-image})
and~$\gamma$ (the {\it complement\/}) be defined by
$$
{\gamma}\,\sigma(i):=n+1-\sigma(i);\quad
{\rho}\,\sigma(i):=\sigma(n+1-i)\quad (1\le i\le n).
$$
The transformation ${\rho\gamma}$ is a bijection
of ${\cal RA}_{2n+1}$ onto ${\cal F\!A}_{2n+1}$
preserving the number of inversions. This makes up a
combinatorial proof of the second identity in (7.7).

\proclaim Proposition 10.1.  We have 
$$\leqalignno{
q^{(2n+1)(2n)/2}A_{2n+1}(q^{-1})&=A_{2n+1}(q);&(7.11)\cr
q^{(2n)(2n-1)/2}A_{2n}(q^{-1})&=A^{\Sec}_{2n}(q).&(7.12)\cr
}
$$

{\it Proof}.
The transformation ${\rho}$ is a bijection
of ${\cal R\!A}_{2n+1}$ onto ${\cal R\!A}_{2n+1}$ with the
property that:
$\inv\sigma+\inv{\rho}\,\sigma=(2n)(2n+1)/2$.
By (7.7) we have
$A_{2n+1}(q)=\sum\limits_{\sigma\in {\cal
F\!A}_{2n+1}}\kern-5ptq^{\inv\sigma} =
\kern-5pt
\sum\limits_{\sigma\in {\cal
F\!A}_{2n+1}}\kern-5pt q^{\inv{\rho}\,\sigma}
=q^{2n(2n+1)/2}\kern-8pt
\sum\limits_{\sigma\in {\cal F\!A}_{2n+1}}\kern-5pt
q^{-\inv\,\sigma} =q^{2n(2n+1)/2}A_{2n+1}(q^{-1})$.
In the same manner, the transformation ${\rho}$ is a bijection
of ${\cal R\!A}_{2n}$ onto ${\cal F\!A}_{2n}$ with the
property that:
$\inv\sigma+\inv{\rho}\,\sigma=(2n)(2n-1)/2$.
Hence, $A_{2n}^{\Sec}(q)=\sum\limits_{\sigma\in {\cal
F\!A}_{2n}}\kern-5ptq^{\inv\sigma} =
\kern-5pt
\sum\limits_{\sigma\in {\cal
R\!A}_{2n}}\kern-5pt q^{\inv{\rho}\,\sigma}
=q^{2n(2n-1)/2}\kern-8pt
\sum\limits_{\sigma\in {\cal R\!A}_{2n}}\kern-5pt
q^{-\inv\,\sigma} =q^{2n(2n-1)/2}A_{2n}(q^{-1})$
\qed

\medskip
\centerline{\bf 8. Proof of Theorem 1.4}

\medskip

By means of the second fundamental transformation
$\Phi$, already mentioned in (7.10), and the bijection
${\bf i}$ that maps each
permutation~$\sigma$ from the symmetric group ${\goth
S}_n$ onto its inverse $\sigma^{-1}$, we can form
$\psi:={\bf i}\,\Phi\,{\bf i}$. The latter bijection has
the following properties:
$$
\Ligne \psi(\sigma)=\Ligne\sigma;
\qquad
\inv\psi(\sigma)=\imaj\sigma.
$$
If ${\bf c}=(c_0,c_1,\ldots,c_m)$ is a composition
of~$n$ and $\sigma=\sigma(1)\sigma(2)\cdots\sigma(n)$ 
a permutation from ${\goth S}_n$, let
$$\leqalignno{
(\Ligne\setminus{\bf c})\sigma
:=\Ligne\sigma&\setminus\{c_0,c_0+c_1,\ldots,
c_0+c_1+\cdots+c_{m-1}\},\cr
\noalign{\hbox{so that for every composition ${\bf c}$ of
$n$ }} 
(\Ligne\setminus{\bf c})
\psi(\sigma)&=(\Ligne\setminus{\bf c})\sigma;
\qquad
\inv\psi(\sigma)=\imaj\sigma.&(8.1)\cr}
$$

Now, consider a $t$-permutation $w=(w_0,w_1,\ldots,
w_m)$ of order~$n$ such that $\Lambda w=
{\bf c}=(c_0,c_1,\ldots,c_m)$. The transformation~$\psi$
maps the {\it permutation} $\sigma:=w_0w_1\cdots w_m$
to another permutation $\psi(\sigma)$.
Let
$(w'_0,w'_1,\ldots,w'_m)$ be the factorization of
$\psi(\sigma)$, written as a word of~$n$ letters, defined
by:
$\lambda w'_0=c_0$, $\lambda w'_1=c_1$, \dots~, 
$\lambda w'_m=c_m$. Property (8.1) implies that the
mapping 
$$
\psi:w=(w_0,w_1,\ldots,w_m)\mapsto w'
:=(w'_0,w'_1,\ldots,w'_m)
$$
is a bijection of ${\cal T}_n$ onto itself having the
properties:
$$
\Lambda w'=\Lambda w={\bf c},\quad
\inv w'=\imaj w.\leqno(8.2)
$$

\goodbreak
\smallskip
For instance, (see [FS78, p. 147] where the same numerical
example is here reproduced) we have:
$\Phi(7\,4\,9\,2\,6\,1\,5\,8\,3)=
4\,7\,2\,6\,1\,9\,5\,8\,3$. Hence,
$$w\!=\!6\,4\,9\,2\,7\,5\,1\,8\,3\buildrel{\bf i}
\over \longrightarrow 
7\,4\,9\,2\,6\,1\,5\,8\,3
\buildrel\Phi
\over \longrightarrow 
4\,7\,2\,6\,1\,9\,5\,8\,3
\buildrel{\bf i}
\over \longrightarrow 
5\,3\,9\,1\,7\,4\,2\,8\,6\!=\!w'
.$$
The $t$-permutation $w=(\epsilon,6, 4, 927, 5183)$ is then
mapped under~$\psi$ onto $w'=(\epsilon, 5, 3, 917, 4286)$.
Moreover, $\imaj w=1+3+5+8=17=1+3+1+3+5+1+3=\inv
w'$.

\goodbreak
\medskip
From Theorem 1.3 and (8.2) it then follows that
$$
A_{n,{\bf c}}(q)=
\sum_{w\in {\cal T}_n,\,
\Lambda w={\bf c}} q^{\inv w}
=\sum_{w\in {\cal T}_n,\,
\Lambda w={\bf c}} q^{\imaj w}.
$$
Now, sum the previous identity over all $t$-compositions
of~$n$. We get:
$$\leqalignno{
\sum_{\mu\,{\bf c}=m}A_{n,{\bf c}}(q)
&=\sum_{\mu\,{\bf c}=m}\;
\sum_{w\in {\cal T}_n,\,
\Lambda w={\bf c}} q^{\inv w}
=\sum_{\mu\,{\bf c}=m}\;
\sum_{w\in {\cal T}_n,\,
\Lambda w={\bf c}} q^{\imaj w}\cr
&=\sum_{\scriptstyle k\ge 0,\,a+b=m
\atop \scriptstyle w\in {\cal T}_{n,k,a,b}}q^{\imaj w}
=\sum_{\scriptstyle k\ge 0\atop \scriptstyle a+b=m}
A_{n,k,a,b}(q).\cr}
$$
This proves (1.17). 

\medskip
For instance, from the tables~2 and 3 we can verify that (1.17)
holds for $n=3$ and
for $m=2$: $A_{3,0,0,2}(q)+A_{3,1,0,2}(q)+A_{3,1,1,1}(q)
+A_{3,1,2,0}(q)+A_{3,2,2,0}(q)
=1+q^2+(2q+2q^2)+q+q^3
=(q+q^2)+(1+q+q^2)+(q+q^2+q^3)
=A_{3,(3)}(q)+A_{3,(2,1,0)}(q)
+A_{3,(0,1,2)}(q)$.

\medskip
The same technique of proof can be used for the
$s$-permutations. We get:
$$\leqalignno{
\sum_{{\bf c}\in\Theta_{n, m+1}^-}A_{n,{\bf c}}(q)
&=\sum_{{\bf c}\in\Theta_{n,m+1}^-}\;
\sum_{w\in {\cal T}_n^-,\,
\Lambda w={\bf c}} q^{\inv w}
=\sum_{{\bf c}\in\Theta_{n,m+1}^-}\;
\sum_{w\in {\cal T}_n^-,\,
\Lambda w={\bf c}} q^{\imaj w}\cr
&=\sum_{\scriptstyle k\ge 0,\,a+b=m
\atop \scriptstyle w\in {\cal T}_{n,k,a,b+1}^-}q^{\imaj w}
=\sum_{\scriptstyle k\ge 0\atop \scriptstyle a+b=m}
B_{n,k,a,b}(q).\cr}
$$
This proves (1.18).

\goodbreak
\centerline{\bf 9. Proof of Theorem 1.5} 
\medskip
If ${\cal J}=(J_0,J_1,\ldots,J_{m-1},J_m)$ is a sequence of
disjoint subsets of the interval
$[\,n\,]:=\{1,2,\ldots,n\}$ of union~$[\,n\,]$ with $m\ge
1$, then $(\#J_0,\#J_1,\ldots,\#J_m)$ is a
composition ${\bf c}\!=\!(c_0,c_1,\ldots, c_m)$
of~$n$. We then write
$\#\,{\cal J}\!:=\!\nobreak {\bf c}$. Also, let
$\inv {\cal J}$ denote the number of ordered pairs $(x,y)$
where
$x\in J_k$, $y\in J_l$, $k<l$ and
$x>y$. A classical result that goes back to MacMahon (see,
e.g. [An76, \S\kern2pt 3.4]) makes it possible to write for
each composition ${\bf c}$ of~$n$
$$\displaylines{\kern3.4cm
\sum_{{\cal J}\!,\,\#{\cal J}={\bf c}}q^{\inv {\cal J}}
={n\brack c_0,c_1,\ldots, c_m}_q,
\hfill\cr
\noalign{\hbox{where the right-hand side is the
$q$-{\it multinomial coefficient} equal to}}
{(q;q)_n\over
(q;q)_{c_0}\,(q;q)_{c_1}\,\cdots\,
\,(q;q)_{c_m}}.\cr
}
$$

\proclaim Theorem 9.1. For each
$t$-composition 
${\bf c}=(c_0,c_1,\ldots,c_{m-1},c_m)$ of~$n$ ($n\geq 1$) we have:
$$
\leqalignno{
A_{n,{\bf c}}(q)&=
{\textstyle{n\brack \,c_0,c_1,\ldots,c_{m-1},
c_m}_q}
A_{c_0}(q)A_{c_1}(q)
\cdots A_{c_{m-1}}(q)
A^{\Sec}_{c_m}(q).\cr
}
$$

{\it Proof}.\quad
Each $t$-permutation $w$ from
${\cal T}_{n}$ such that $\Lambda w={\bf c}$ and
$\mu\,{\bf c}=m\ge 1$ is completely characterized by a
sequence
$$\bigl(\,(I_0,\sigma_0), (I_1,\sigma_1),
\ldots,(I_m,\sigma_m)\,\bigr),$$ 
having the following properties:

(i) the sequence
${\cal I}(w):=(I_0,I_1,\ldots,I_m)$ consists of
disjoint subsets of the interval
$[\,n\,]:=\{1,2,\ldots,n\}$ of union~$[\,n\,]$; moreover,
$\#{\cal I}(w):={\bf c}\,$;

(ii) $\sigma_0\in {\cal RA}_{c_0}$,
$\sigma_1\in {\cal F\!A}_{c_1}$,
\dots~,  
$\sigma_m\in {\cal F\!A}_{c_m}$.

\medskip
\noindent
 If ${\cal I}(w)={\cal J}$, then
$\inv w=\inv {\cal J}+\inv\sigma_0
+\inv\sigma_1+\cdots+\inv\sigma_{m-1}+\inv\sigma_m$.
Hence,
$$\leqalignno{
A_{n,{\bf c}}(q)&=\sum_{w\in {\cal T}_n,\Lambda w={\bf c}}
q^{\inv w}=
\sum_{w\in {\cal T}_n,\,\#{\cal I}(w)={\bf c}}\kern-15pt
q^{\inv w}\cr
& =
\sum_{{\cal J}\!,\,\#{\cal J}={\bf c}}\;\sum_{w,\,{\cal
I}(w)={\cal J}} q^{\inv w}\cr
&=
\sum_{{\cal J}\!,\,\#{\cal J}={\bf c}}q^{\inv {\cal J}}
\kern-8pt 
\sum_{\sigma_0\in {\cal RA}_{c_0}}q^{\inv \sigma_0}
\kern-8pt 
\sum_{\sigma_1\in {\cal F\!A}_{c_1}}
q^{\inv\sigma_1}
\times\cdots
\times\kern-10pt
\sum_{\sigma_m\in {\cal F\!A}_{c_m}}
\kern-10pt 
q^{\inv\sigma_m}\cr
&=
{\textstyle{n\brack \,c_0,c_1,\ldots,
c_m}_q}
A_{c_0}(q)A_{c_1}(q)
\cdots
 A^{\Sec}_{c_m}(q).\qed\cr}
$$

\medskip
The factorial generating functions for the polynomials
$$\leqalignno{
A_n(x,q)&:=\sum_{{\bf c}\in \Theta_n}x^{\mu\,{\bf
c}}A_{n,{\bf c}}(q),\qquad
B_n(x,q):=\sum_{{\bf c}\in \Theta_n^-}x^{\mu\,{\bf
c}-1}B_{n,{\bf c}}(q),
\cr}
$$
can be derived from Theorem 9.1. 

\medskip

{\it Proof of Theorem 1.5}.\quad
For $m\ge 1$ we  have:
$$\displaylines{\quad \sum_{n\ge 0}{u^n\over (q;q)_n}
\sum_{\scriptstyle w\in \Theta_n\atop
\scriptstyle \mu\,{\bf c}=m}A_{n,{\bf c}}(q)
\hfill\cr
\kern1cm{}=\sum_{n\ge 0}{u^n\over (q;q)_n} 
\sum_{\scriptstyle c_0+c_1+\cdots\ \atop
\scriptstyle \ +c_m=n}\kern-12pt  
{\textstyle{n\brack \,c_0,c_1,\ldots,
c_m}_q}
A_{c_0}(q)A_{c_1}(q)
\cdots
 A^{\Sec}_{c_m}(q)
\hfill \cr
\kern1cm{}=\sum_{n\ge 0}u^n
\sum_{\scriptstyle c_0+c_1+\cdots\ \atop
\scriptstyle \ +c_m=n} {A_{c_0}(q)\over
(q;q)_{c_0}}{A_{c_1}(q)\over (q;q)_{c_1}}\cdots
{A_{c_{m-1}}(q)\over
(q;q)_{c_{m-1}+1}}{A^{\Sec}_{c_m}(q)\over (q;q)_{c_m}}.
\hfill\cr }
$$
When $m\ge 1$, the integers $c_0$ and $c_m$ are even; if,
furthermore, $m\ge 2$, then $c_1$, \dots,~$c_{m-1}$ are
odd. Hence,
the previous identity may be rewritten as
$$\displaylines{\sum_{n\ge 0}{u^n\over (q;q)_n}
\sum_{\scriptstyle w\in \Theta_n\atop
\scriptstyle \mu\,{\bf c}=m}A_{n,{\bf c}}(q) \hfill\cr
\quad{}
=
\sum_{j_0\ge 0}A_{2j_0}(q){u^{2j_0}\over (q;q)_{2j_0}}\;
\Bigl(\sum_{j\ge 0}A_{2j+1}(q){u^{2j+1}\over
(q;q)_{2j +1}}\Bigr)^{m-1}
\sum_{j_m\ge 0}A^{\Sec}_{2j_m}(q){u^{2j_m}\over
(q;q)_{2j_m}}\hfill\cr
\quad{}
=\sec_q(u)\, \bigl(\tan_q(u)\bigr)^{m-1}\,
\Sec_q(u).\hfill\cr }
$$
When $m=0$, then $n$ is {\it odd}. Hence, $A_{n,{\bf c}}(q)
=A_n(q)$ ($n$ odd) and
$$
\sum_{n\ge 0}{u^n\over (q;q)_n}
\sum_{\scriptstyle w\in \Theta_n\atop
\scriptstyle \mu\,{\bf c}=0}A_{n,{\bf c}}(q) 
=\tan_q(u).
$$
Thus,
$$
\eqalignno{\sum_{n\ge 0}A_n(x,q){u^n\over (q;q)_n}
&=\sum_{n\ge 0}{u^n\over (q;q)_n}
\sum_{m\ge 0}x^m\sum_{\scriptstyle w\in \Theta_n\atop
\scriptstyle \mu\,{\bf c}=m}A_{n,{\bf c}}(q) 
\cr
&=
\sum_{n\ge 0}{u^n\over (q;q)_n}\kern-3pt
\sum_{\scriptstyle w\in \Theta_n\atop
\scriptstyle \mu\,{\bf c}=0}\kern-5pt A_{n,{\bf c}}(q)
+\kern-3pt 
\sum_{m\ge 1}x^m\sum_{n\ge 0}{u^n\over (q;q)_n}
\kern-3pt 
\sum_{\scriptstyle w\in \Theta_n\atop
\scriptstyle \mu\,{\bf c}=m}\kern-5pt A_{n,{\bf c}}(q) \cr
&=\tan_q(u)+\sum_{m\ge 1}\sec_q(u)\,
\bigl(x\tan_q(u)\bigr)^{m-1}\,
x\Sec_q(u)\cr
&=\tan_q(u)+\sec_q(u)(1-x\tan_q(u))^{-1}x\Sec_q(u),\cr}
$$
which proves (1.19).
\goodbreak

The proof of (1.20) is quite similar. The only difference is
the fact that ${\bf c}=(c_0,c_1,\ldots,c_m,0)$ is now an
$s$-composition, so that, $c_0$ is even and if $m\ge 1$
all the other~$c_i$'s are odd. We then get:
$$
\leqalignno{\sum_{n\ge 0}{u^n\over (q;q)_n}
\sum_{\scriptstyle w\in \Theta_n^-\atop
\scriptstyle \mu\,{\bf c}=m+1}A_{n,{\bf c}}(q)
&=
\sum_{j_0\ge 0}A_{2j_0}(q){u^{2j_0}\over (q;q)_{2j_0}}\;
\Bigl(\sum_{j\ge 0}A_{2j+1}(q){u^{2j+1}\over
(q;q)_{2j +1}}\Bigr)^{m}\cr
&=\sec_q(u)\, \bigl(\tan_q(u)\bigr)^{m}.\cr 
\noalign{\hbox{Then}}
\sum_{n\ge 0}B_n(x,q){u^n\over (q;q)_n}
&=
\sum_{m\ge 0}x^m\sum_{n\ge 0}{u^n\over (q;q)_n}
\kern-5pt 
\sum_{\scriptstyle w\in \Theta_n^-\atop
\scriptstyle \mu\,{\bf c}=m+1}A_{n,{\bf c}}(q)\cr
&=\sum_{m\ge 0}\sec_q(u)\,
\bigl(x\tan_q(u)\bigr)^{m}\cr
&=\sec_q(u)(1-x\tan_q(u))^{-1},\cr}
$$
which proves (1.20).\qed

\bigskip
\goodbreak
\centerline{\bf 10. Specializations} 

\medskip
Our three families of polynomials 
$(A_{n,k,a,b}(q))$, $(B_{n,k,a,b}(q))$,
and $(A_{n,{\bf c}}(q))$
involve specializations that relate to other classes of
generating polynomials or classical numbers that have been
studied in previous works. Those polynomials are
displayed in Tables 2--4 at the end of the
paper. Each table appears as a matrix, whose
$(n,m)$-cell contains several polynomials. In the
$(n,m)$-cell of Table~2 (resp. Table~3) are reproduced all
the polynomials 
$A_{n,k,a,b}(q)$ (or $B_{n,k,a,b}(q)$) such that $a+b=m$
(resp. $A_{n,{\bf c}}(q)$ (or $B_{n,{\bf c}}(q)$) such that 
$\mu\,{\bf c}=m$.
The specializations we deal with
refer to rows, columns or diagonals of those tables. Others are
obtained by summing the above polynomials with respect to
certain subscripts.

To this end we use the following notations:
$$\leqalignno{
A_{n,k,a+b=m}(q)&:=\sum_{a+b=m}A_{n,k,a,b}(q);&(10.1)\cr
A_{n,a+b=m}(q)&:=\sum_{\scriptstyle k\ge 0\atop
\scriptstyle a+b=m}A_{n,k,a,b}(q);&(10.2)\cr
\noalign{\hbox{By Theorem 1.4 we also have:}}
A_{n,a+b=m}(q)&:=\sum_{\mu\,{\bf c}=m}A_{n,{\bf c}}(q).\cr}
$$
Analogous definitions are made for the polynomials
$B_{n,k,a,b}(q)$.

\medskip
10.1. {\it The first column of $(A_{n,k,a,b}(q))$,
$(B_{n,k,a,b}(q))$}.\quad
The $(t,q)$-analogs of tangent and secant have been
introduced in our previous paper  [FH11]. For each $r\ge 0$
form the
$q$-series:
$$
\leqalignno{
\qquad\quad\sin_q^{(r)}(u)&:=
\sum_{n\ge 0}(-1)^n
{(q^r;q)_{2n+1}\over (q;q)_{2n+1}}u^{2n+1};\cr
\cos_q^{(r)}(u)&:=
\sum_{n\ge 0}(-1)^n
{(q^r;q)_{2n}\over (q;q)_{2n}}u^{2n};\cr
\tan_q^{(r)}(u)&:={\sin_q^{(r)}(u)\over
\cos_q^{(r)}(u)};\cr
\sec_q^{(r)}(u)&:={1\over \cos_q^{(r)}(u)}.\cr
}
$$
The $(t,q)$-{\it analogs of the tangent and
secant numbers} have been defined as the coefficients
$T_{2n+1}(t,q)$ and
$E_{2n}(t,q)$, respectively, in the following two series:
$$\leqalignno{
\sum_{r\ge 0}t^r\tan_q^{(r)}(u)
&=\sum_{n\ge 0}{u^{2n+1}\over
(t;q)_{2n+2}}T_{2n+1}(t,q) ;&(10.3)\cr
\sum_{r\ge 0}t^r\sec_q^{(r)}(u)&=\sum_{n\ge
0}{u^{2n}\over (t;q)_{2n+1}}E_{2n}(t,q).&(10.4)\cr
}
$$
It was then proved that 
$T_{2n+1}(t,q)$ and $E_{2n}(t,q)$  have the following
combinatorial interpretations:
$$\leqalignno{
T_{2n+1}(t,q)&=\sum_{\sigma\in {\cal
RA}_{2n+1}}t^{1+\ides\sigma}q^{\imaj\sigma};\cr
E_{2n}(t,q)&=\sum_{\sigma\in {\cal
RA}_{2n}}t^{1+\ides\sigma}q^{\imaj\sigma}.\cr}
$$

Now, the set ${\cal T}_{2n+1,k,0,0}$ is the set of all
$t$-permutations~$w=(w_0)$ of order $(2n+1)$, where
$w_0$ is simply an element of ${\cal RA}_{2n+1}$, that is,
a rising alternating permutation of order~$(2n+1)$. Hence,
$$
T_{2n+1}(t,q)=\sum_{k\ge 1}t^{k+1}A_{2n+1,k,0,0}(q).
$$
Accordingly, Theorem 4.1 provides a method for calculating the 
polynomials $T_{2n+1}(t,q)$, only defined so far by their
generating function (10.3). In an equivalent manner, we can
also say that the factorial generating function for the first
column of the matrix $(A_{n,k,a,b}(q))$ is given by
$$
\sum_{r\ge 0}t^r\tan_q^{(r)}(u)
=\sum_{n\ge 0}{u^{2n+1}\over
(t;q)_{2n+2}}\sum_{k\ge 1}t^{k+1}A_{2n+1,k,0,0}(q).
$$
In the same way, we get:
$$
E_{2n}(t,q)=\sum_{k\ge 1}t^{k+1}B_{2n,k,0,0}(q),
$$
which also provides, either a way of calculating the
$(t,q)$-analogs
$E_{2n}(t,q)$ of the secant numbers, or writing the
factorial generating function for the first column of the
matrix $(B_{n,k,a,b}(q))$.

\goodbreak
\medskip
10.2. {\it The super-diagonal $(n,n+1)$ of the matrix 
$(A_{n,k,a,b}(q))$}.\quad As will be shown, the polynomials
$A_{n,k,a,b}(q)$ $(a+b=n+1)$ of that super-diagonal
provide a {\it refinement} of the {\it Carlitz
$q$-analogs of the Eulerian polynomials} [Ca54].
Twenty-one years later [Ca75] Carlitz also showed that they
were generating polynomials for the symmetric groups by
the pair ``des" (number of descents) and ``maj" (major
index). Let
$(A_n(t,q))$ be the sequence of those polynomials, written as
$A_n(t,q)=\smash{\sum\limits_{j\ge 0}}A_{n,j}(q)$ $(n\ge
0)$.  The recurrence
$$
A_{n,j}(q)
=(1\!\!+\!q\!+\!\cdots
\!+\!\!q^j)\,A_{n-1,j}(q)\!+\!(q^j\!\!
+\!q^{j\!+\!1}
\!\!+\!\cdots\!+\!q^{n-1})
\,A_{n-1,j-1}(q),
$$
with the initial conditions
$A_{0,j}(q)=A_{1,j}(q)=\delta_{0,j}$, provides a method for
calculating them.

Their first values 
are reproduced in the following table:
\medskip
\noindent
$A_0(t,q)=A_1(t,q)=1$; 
$A_2(t,q)=1+tq$; 
$A_3(t,q)=1+2tq(q+1)+t^2q^3$;

\noindent
$A_4(t,q)=1+tq(3q^2+5q+3)+t^2q^3(3q^2+5q+3)+t^3q^6$;
\hfil\break
\noindent
$A_5(t,q)=1+tq(4q^3+9q^2+9q+4)+t^2q^3(6q^4+16q^3+22q^2+16q+6)
+\null\hskip 1.6cm t^3q^6(4q^3+9q^2+9q+4)+t^4q^{10}$.

\medskip
Now, go back to 
the recurrence for the polynmials
$A_{n,k,a,b}(q)$ shown in (4.1) and rewrite it
when $a'+b'=m'=n+2$. The coefficients
$A_{n,k'-1,a,m'+1-a}(q)
=A_{n,k'-1,a,n+3-a}(q)$ and
$A_{n,k',a,m'+1-a}(q)=A_{n,k',a,n+3-a}(q)$ vanish,
because $A_{n,k,a,b}(q)=0$ when $a+b\ge n+2$. Hence,
$$
\displaylines{ 
A_{n+1,k',a',n+2-a'}(q)\!=\!q^{k'}\Bigl(\,
\sum_{0\le a\le a'-1\le n}\kern-10pt 
A_{n,k'-1,a,n+1-a}(q)\hfill\cr
\noalign{\vskip-10pt}
\hfill{}+ \sum_{1\le a'\le a\le n+1}\kern-10pt 
A_{n,k',a,n+1-a}(q)\Bigr),\quad\cr
}
$$
valid for $n\ge 0$ with the initial condition:
$A_{0,k,a,b}(q)=\delta_{k,0}\delta_{a,1}\delta_{b,0}$.
For $n\ge 0$ let $A_{n,k,a}(q):=A_{n,k,a,n+1-a}(q)$, so that
the previous recurrence can be written in the form 
$$
A_{n+1,k',a'}(q)\!=\!q^{k'}\Bigl(\,
\sum_{0\le a\le a'-1\le n}\kern-10pt 
A_{n,k'-1,a}(q)+ \sum_{1\le a'\le a\le n+1}\kern-10pt 
A_{n,k',a}(q)\Bigr),
$$
with $0\le a'\le n+2$ and makes it possible the
calculations of those polynomials:

\medskip
$A_{1,0,1}(q)=1$;

\smallskip
$A_{2,0,1}(q)=1$, $A_{2,1,2}(q)=q$;

\smallskip
$A_{3,0,1}(q)\!=\!1$, $A_{3,1,1}(q)\!=\!q^2$,
$A_{3,1,2}(q)\!=\!q+q^2$, $A_{3,1,3}(q)\!=\!q$, 
$A_{3,2,3}(q)\!=\!q^3$;

\smallskip
$A_{4,0,1}(q)=1$, $A_{4,1,1}(q)=2q^2+2q^3$,
$A_{4,1,2}(q)=q+2q^2+q^3$, $A_{4,1,3}(q)=q+q^2$,

$A_{4,1,4}(q)=q$, $A_{4,2,1}(q)=q^5$,
$A_{4,2,2}(q)=q^4+q^5$,
$A_{4,2,3}(q)=q^3+2q^4+q^5$,\hfil\break\indent
$A_{4,2,4}(q)=2q^3+2q^4$, $A_{4,3,4}(q)=q^6$.

\proclaim Theorem 10.1. For each $n\ge 0$ 
and each $j\ge 0$ the coefficient $A_{n,j}(q)$ of~$t^j$ in
the Carlitz $q$-Eulerian polynomial $A_n(t,q)$ admits the
following refinement:
$$A_{n,j}(q)=\sum_{a\ge 0}A_{n,j,a}(q).
\leqno(10.5)
$$

For instance,
$A_{4,1}(q)=3q\!+\!5q^2\!+\!3q^3=(2q^2\!+\!2q^3)\!+\!(q\!+\!2q^2\!+\!q^3)
\!+\!(q\!+\!q^2)\!+\!q=A_{4,1,1}(q)\!+\!A_{4,1,2}(q)\!+\!A_{4,1,3}(q)
+\!A_{4,1,4}(q)$.

\medskip
{\it Proof}.\quad
Each polynomial $A_n(t,q)$ is the generating polynomial for
${\goth S}_n$ by the pair (``number of descents",``major
index"), as established by Carlitz [Ca75],~or, in an
equivalent manner, by the pair
$(\ides,\imaj)$. This can also be expressed~by
$$
A_{n,j}(q)=\sum_{\sigma\in {\goth S}_n,\ides\sigma=j}
q^{\imaj\sigma},\leqno(10.6)
$$ 

By Theorem~9.1 
$A_{n,j,a}(q)=A_{n,j,a,n+1-a}(q)$ is the
generating polynomial for the set of all
$t$-permutations~$w$ of order~$n$ such that $\mu\,
w=n+1$, $\ides w=j$ and $\min w=a$ by ``imaj."
Such  $t$-permutations are of the form
$w=(w_0,w_1,w_2,\ldots, w_{n+1})$, so that necessarily,
$w_0=w_{n+1}=\epsilon$, and the other components $w_i$
are one-letter words. Accordingly, $A_{n,j,a}(q)$ is the
generating polynomial for all (ordinary) {\it permutations}
$\sigma=w_1w_2\cdots w_n$ of $12\cdots n$ such that 
$\ides\sigma=j$ and $w_a=1$, that is,
$$
A_{n,j,a}(q)=\sum_{\scriptstyle \sigma\in {\goth
S}_n,\ides\sigma=j,\atop\scriptstyle \min \sigma=a}
q^{\imaj\sigma}.\leqno(10.7)
$$
Thus (10.5) is a consequence of (10.6) and (10.7).\qed

\medskip
For each integer $n\ge 1$ let $[n]_q:=(q;q)_n/(1-q)^n
=1+q+\cdots+q^{n-1}$. By summing the $A_{n,j,a}$'s over
the pair $(j,a)$ we get the polynomial
$A_{n,a+b=n+1}(q)$ defined in (10.2), which is the
generating polynomial for ${\goth S}_n$ by ``inv,"
well-known to be equal to
$$
A_{n,a+b=n+1}(q)=[1]_q\,[2]_q\,\cdots\,[n]_q,
\leqno(10.8)
$$
also equal (using the same combinatorial interpretation) to
$B_{n,a+b=n}(q)=\sum\limits_{k\ge 0,a+b=n}B_{n,k,a,b}(q)$.

\bigskip
10.3. {\it The subdiagonal $(n,n-1)$ of the matrix
$(A_{n,k,a,b}(q))$}.\quad
Our purpose is to evaluate the polynomial
$A_{n,a+b=n-1}(q)$ for each $n\ge 1$, which is the
generating function for all $t$-permutations of
order~$n$ such that $\mu\,w=n-1$ by ``inv." Such
$t$-permutations~$w$ have one of the three forms:

(1) $w=(x_1x_2,x_3,x_4,\ldots,x_n,\epsilon)$ with
$x_1<x_2$;

(2) $w=(\epsilon,x_1,x_2,\ldots,x_{n-2},x_{n-1}x_n)$ with
$x_{n-1}>x_n$;

(3) $(\epsilon,x_1,\ldots,x_{i-1},x_ix_{i+1}x_{i+2},
x_{i+3},\ldots,x_n,\epsilon)$ with $1\le i\le n-2$,
$x_i>x_{i+1}$, $x_{i+1}<x_{i+2}$.

\medskip
The g.f. of the $t$-permutations of form~(1) or~(2) by
``inv" is equal to
$$
(1+q){n\brack
2,1^{n-2}}_q=[2]_q\,[3]_q\,[4]_q\,\cdots\,[n]_q.
$$

The g.f. of the $t$-permutations of form~(3) by ``inv" is
equal to
$$\leqalignno{
\sum_{1\le i\le n-2}
{n\brack 1^{i-1},3,1^{n-i-2}}_q\times q(q+1)&=
(n-2)q(q+1){n\brack 3,1^{n-3}}_q\cr
&=(n-2)q\,[2]_q\,[4]_q\,\cdots\,[n]_q,\cr}
$$
so that the total g.f. is equal to
$$\displaylines{\quad
[2]_q\,[4]_q\,\cdots\,[n]_q([3]_q+(n-2)q)\hfill\cr
\kern2cm{}=[2]_q\,[4]_q\,\cdots\,[n]_q(1+q+q^2+(n-2)q)
\hfill\cr
\kern2cm{}=[1]_q\,[2]_q\bigl(1+(n-1)q+q^2\bigr)\,
[4]_q\,\cdots\,[n]_q.\hfill\cr
}$$
Thus,
$A_{n,a+b=n-1}(q)= [1]_q\,[2]_q\bigl(1+(n-1)q+q^2\bigr)\,
[4]_q\,\cdots\,[n]_q $.

\medskip
10.4. {\it The subdiagonal $(n,n-2)$ of the matrix
$(B_{n,k,a,b}(q))$}.\quad Using the same combinatorial
technique as in 10.3, but this time operating with the
$s$-permutations we get the following evaluation for each
$n\ge 2$:
$$
B_{n,a+b=n-2}(q)=
(1+(n-1)q+(n-1)q^2)\,[4]_q\,\cdots\,[n]_q.
$$

\medskip
10.5. {\it Two $q$-analogs of the Springer numbers}.\quad
It was recalled in the Introduction ((1.7) and (1.20)) that
$\sec(u)(1-x\tan (u))^{-1}$ for $x=1$ was the exponential
generating function for the {\it Springer numbers}.
Referring to (1.21) we then see that
$\sec_q(u)(1-\tan_q(u))^{-1}$ is the factorial generating
function for the $q$-analogs of the Springer numbers,
which are simply the generating polynomials for the
$s$-permutations by ``imaj" or ``inv."

Note that $\Sec_q(u)(1-\tan_q(u))^{-1}$ is also the
factorial generating for such $q$-analogs, but this time
for the $S$-permutations by ``imaj" or ``inv;"

\bigskip
10.6. {\it $t$-compositions and Fibonacci triangle}.\quad
Let $\alpha(n,m):=\#\Theta_{n,m}$ (resp. $\beta(n,m):=\#\Theta_{n,m+1}^-$) 
	be the {\it
number} of  
$t$-compositions (resp. $s$-compositons).
From the previous lists of the $\Theta_i$'s 
made in Section~1 we have the next
table, where the 
$\alpha(n,m)$'s (resp. $\beta(n,m)$'s) have been reproduced
in bold face (resp. plain type). To the right are displayed
the row sums of those entries, which will be proven to be the classical 
Fibonacci numbers.
$$
\bordermatrix{m=& 0&1&2&3&4&5&6&7\cr
n=0&1&\bf1\cr
\qquad 1&\bf 1&1&\bf1\cr
\qquad 2&1&\bf2&1&\bf1\cr
\qquad 3&\bf1&2&\bf3&1&\bf1\cr
\qquad 4&1&\bf3&3&\bf4&1&\bf1\cr
\qquad 5&\bf1&3&\bf6&4&\bf5&1&\bf1\cr
\qquad
6&1&\bf4&6&\bf10&5&\bf6&1&\bf1\cr}
\matrix{\bf 1&1\cr
\bf 2&1\cr
\bf 3&2\cr
\bf 5&3\cr
\bf 8&5\cr
\bf 13&8\cr
\bf 21&13\cr 
}
$$
\vskip-2pt
\centerline{Fig. 10.1. The coefficients $\bf \alpha(n,m)$ and
$\beta(n,m)$}
\goodbreak
\medskip

The mapping 
$ (c_0,c_1,\ldots,c_m, 0) \mapsto (c_0,c_1,\ldots,c_m-1) $
is a bijection of $ \Theta_{n,m+1}^-$ onto $\Theta_{n-1,m}$,
because $c_m$ is odd. Hence
$\beta(n,m)=\alpha(n-1,m)\quad(n\ge 1)$.
We now only study the numbers $\alpha(n,m)$.

\proclaim Proposition 10.2. With the initial values
$\alpha(0,m)=\delta_{1,m}$,
$\alpha(1,m)=\delta_{0,m}+\delta_{2,m}$, the entries
$\alpha(n,m)$ are inductively given by
$$
\alpha(n,m)=\alpha(n-1,m-1)+\alpha(n-2,m)
\quad(n\ge 2).
$$

{\it Proof}.\quad
Let ${\bf c}=(c_0, c_1, \ldots, c_{m-1}, c_m)\in \Theta_{n,m}$.  
If $c_m=0$ we define $\phi({\bf c})=(c_0, c_1, \ldots, c_{m-1}-1) \in\Theta_{n-1, m-1}$. 
If $c_m\geq 1$, then $c_m\geq 2$ because $c_m$ is even. We define $\phi({\bf c})=(c_0, c_1, \ldots, c_{m-1}, c_m-2) \in\Theta_{n-2, m}$. We verify that $\phi$ is a bijection between
$\Theta_{n,m}$ and $\Theta_{n-1, m-1} + \Theta_{n-2, m}$.\qed

\medskip
For each $n\ge 0$ let $A_n(x)=\sum_{m\ge
0}\alpha(n,m)x^m$ be the generating polynomials of the
coefficients $\alpha(n,m)$. From Proposition~10.23 it follows
that $A_0(x)=x$, $A_1(x)=1+x^2$ and the recurrence
formula
$$
A_{n+1}(x)=xA_n(x)+A_{n-1}(x)\quad (n\geq 1).
$$

Let $A(x;u):=\sum_{n\ge 0}A_n(x)u^n$. Then,
$A(x;u)-x-u(1+x^2)=xu(A(x;u)-x)+u^2A(x;u)$, so that
$A(x;u)(1-xu-u^2)=x+u$ and
$$
A(x;u)=\sum_{n\ge 0}A_n(x)u^n
=\sum_{{\bf c}\in \Theta}x^{\mu\,{\bf
c}}u^{|{\bf c}|}={x+u\over 1-u(x+u)}.
$$
Let $x=1$ we obtain the generating function for the row sums
$A_n(1)=\sum \alpha(n,m)$ ($m\geq 0$), which is equal to $(1+u)/(1-u-u^2)$.
Thus, the row sums are the classical Fibonacci numbers.

The polynomials $A_n(x)$ are related to the polynomials $F_n(x)$ already introduced in Sloane's Integer Encyclopedia [Sl06] under reference A102426
by $A_n(x)=x^{n+1}F_{n+1}(x^{-2})$ for $n\ge 0$. Accordingly, the $t$-compositions provide a natural combinatorial interpretation for their coefficients.

\medskip
10.7. {\it Further comment}.\quad   Dominique Dumont [Du12] has drawn our attention
to the two papers by Carlitz-Scoville [CS72] and Fran\c con
[Fr78]. Instead of $t$-permutations or snakes, Carlitz and Scoville
have dealt with ``up-down sequences of length
$n+m$ with $m$ infinite elements." Such a sequence is a rising
alternating permutation $x_1x_2 \cdots x_{n+m}$ containing
all the integers $1,2, \ldots, n$ and $m$ letters equal to $-\infty$. Note that
replacing all the commas in each $t$-permutation $w=(w_0,w_1, \ldots , w_m)$
by  $-\infty$ makes up a bijection of the set of all $t$-permutations onto
the set of all Carlitz-Scoville sequences. In Fran\c con [Fr78]
can be found an unexpected combinatorial interpretation
of the entries $b(n,m)$ in terms of computer file histories.

\bigskip

\centerline{\bf 11. Tables} 

\medskip
Four tables are being displayed, 
the first one 
containing the values of $a(n,m)$
and $b(n,m)$ for $0\le n\le 6$, 		 
		the second one 
containing the values of the polynomials $A_{n,k,a,b}(q)$
and $B_{n,k,a,b}(q)$ for $0\le n\le 4$, the third one 
for the polynomials $A_{n,{\bf c}}(q)$ 
 for $0\le n\le 4$. The last one 
contains the values of the polynomials $A_{n,a+b=n}(q)$
and $B_{n,a+b=n}(q)$, whose definitions are given in
(10.1)--(10.2).

$$
\bordermatrix{m=& 0&1&2&3&4&5&6&7\cr
n=0&1&\bf1\cr
\qquad 1&\bf 1&1&\bf1\cr
\qquad 2&1&\bf2&2&\bf2\cr
\qquad 3&\bf2&5&\bf8&6&\bf6\cr
\qquad 4&5&\bf16&28&\bf40&24&\bf24\cr
\qquad 5&\bf16&61&\bf136&180&\bf240&120&\bf120\cr
\qquad
6&61&\bf272&662&\bf1232&1320&\bf1680&720&\bf720\cr}
\matrix{\bf 1.2^0&1\cr
\bf 1.2^1&1\cr
\bf 1.2^2&3\cr
\bf 2.2^3&11\cr
\bf 5.2^4&57\cr
\bf 16.2^5&361\cr
\bf 61. 2^6&2763\cr 
}
$$
\nobreak
\centerline{Table 1. The coefficients $\bf a(n,m)$ and
$b(n,m)$}

\vfill\eject

{\eightpoint
\centerline{\vtop{\offinterlineskip
\halign{\vrule\hbox{\vrule height 9pt depth 4pt width0pt}
\hfil$#$\hfil\thinspace \vrule 
&\thinspace \hfil$#$\hfil\thinspace \vrule
&\thinspace \hfil$#$\hfil\thinspace \vrule
&\thinspace \hfil$#$\hfil\thinspace \vrule
&\thinspace \hfil$#$\hfil\thinspace \vrule
\cr
\noalign{\hrule}
m\!=\!&0&1&2&3\cr
\noalign{\hrule}
\kern-3pt
n\!=\!0\quad&B_{0,-\!1,0,0}\!\!=\!\!1&
A_{0,0,1,0}\!\!=\!\!1&&\cr
\noalign{\hrule}
\quad1\quad&A_{1,0,0,0}\!\!=\!\!1
&B_{1,0,1,0}\!\!=\!\!1&A_{1,0,1,1}\!\!=\!\!1&\cr
\noalign{\hrule}
2&B_{2,0,0,0}\!\!=\!\!1&A_{2,0,0,1}\!\!=\!\!1&
B_{2,0,1,1}\!\!=\!\!1&A_{2,0,1,2}\!\!=\!\!1\cr
&&A_{2,1,1,0}\!\!=\!\!q&B_{2,1,2,0}\!\!=\!\!q
&A_{2,1,2,1}\!\!=\!\!q\cr
\noalign{\hrule}
3&A_{3,1,0,0}\!\!=\!\!q\!\!+\!\!q^2
&B_{3,0,0,1}\!\!=\!\!1&A_{3,0,0,2}\!\!=\!\!1&B_{3,0,1,2}\!\!=\!\!1\cr
&&B_{3,1,0,1}\!\!=\!\!q^2&A_{3,1,0,2}\!\!=\!\!q^2
&B_{3,1,1,2}\!\!=\!\!q^2\cr
&&B_{3,1,1,0}\!\!=\!\!2q\!\!+\!\!q^2
&A_{3,1,1,1}\!\!=\!\!2q\!\!+\!\!2q^2
&B_{3,1,2,1}\!\!=\!\!q\!\!+\!\!q^2\cr
&&&A_{3,1,2,0}\!\!=\!\!q&B_{3,1,3,0}\!\!=\!\!q\cr
&&&A_{3,2,2,0}\!\!=\!\!q^3&B_{3,2,3,0}\!\!=\!\!q^3\cr
\noalign{\hrule}
4&B_{4,1,0,0}
\!\!=\!\!q\!\!+\!\!2q^2\!\!+\!\!q^3&A_{4,1,0,1}
\!\!=\!\!q\!\!+\!\!3q^2\!\!+\!\!2q^3
&B_{4,0,0,2}\!\!=\!\!1&A_{4,0,0,3}\!\!=\!\!1\cr
&B_{4,2,0,0}\!\!=\!\!q^4&A_{4,1,1,0}
\!\!=\!\!q\!\!+\!\!q^2
&B_{4,1,0,2}\!\!=\!\!2q^2\!\!+\!\!2q^3
&A_{4,1,0,3}\!\!=\!\!2q^2\!\!+\!\!2q^3
\cr
&&A_{4,2,0,1}\!\!=\!\!q^4\!\!+\!\!q^5
&B_{4,1,1,1}\!\!=\!\!2q\!\!+\!\!4q^2\!\!+\!\!2q^3
&A_{4,1,1,2}\!\!=\!\!2q\!\!+\!\!5q^2\!\!+\!\!3q^3
\cr
&&A_{4,2,1,0}\!\!=\!\!2q^3\!\!+\!\!3q^4\!\!+\!\!q^5
&B_{4,1,2,0}\!\!=\!\!2q\!\!+\!\!q^2
&A_{4,1,2,1}\!\!=\!\!2q\!\!+\!\!2q^2
\cr
&&&B_{4,2,0,2}\!\!=\!\!q^5&A_{4,1,3,0}\!\!=\!\!q\cr
&&&B_{4,2,1,1}\!\!=\!\!2q^4\!\!+\!\!q^5
&A_{4,2,0,3}\!\!=q^5\cr
&&&B_{4,2,2,0}\!\!=\!\!3q^3\!\!+\!\!4q^4\!\!+\!\!q^5
&A_{4,2,1,2}\!\!=\!\!2q^4\!\!+\!\!2q^5\cr
&&&&A_{4,2,2,1}\!\!=\!\!3q^3\!\!+\!\!5q^4\!\!+\!\!2q^5\cr
&&&&A_{4,2,3,0}\!\!=\!\!2q^3\!\!+\!\!2q^4\cr
&&&&A_{4,3,3,0}\!\!=\!\!q^6\cr
\noalign{\hrule}
\noalign{\smallskip}
}}}

\medskip
\centerline{\vtop{\offinterlineskip
\halign{\vrule\hbox{\vrule height 9pt depth 4pt width0pt}
\hfil$#$\hfil\thinspace \vrule 
&\thinspace \hfil$#$\hfil\thinspace \vrule
&\thinspace \hfil$#$\hfil\thinspace \vrule
\cr
\noalign{\hrule}
m\!=\!&4&5\cr
\noalign{\hrule}
n\!=\!3\quad
&A_{3,0,1,3}\!\!=\!\!1&\cr
&A_{3,1,1,3}\!\!=\!\!q^2&\cr
&A_{3,1,2,2}\!\!=\!\!q\!\!+\!\!q^2&\cr
&A_{3,1,3,1}\!\!=\!\!q&\cr
&A_{3,2,3,1}\!\!=\!\!q^3&\cr
\noalign{\hrule}
\ 4&B_{4,0,1,3}\!\!=\!\!1&A_{4,0,1,4}\!\!=\!\!1\cr
&B_{4,1,1,3}\!\!=\!\!2q^2\!\!+\!\!2q^3
&A_{4,1,1,4}\!\!=\!\!2q^2\!\!+\!\!2q^3\cr
&B_{4,1,2,2}\!\!=\!\!q\!\!+\!\!2q^2\!\!+\!\!q^3
&A_{4,1,2,3}\!\!=\!\!q\!\!+\!\!2q^2\!\!+\!\!q^3\cr
&B_{4,1,3,1}\!\!=\!\!q\!\!+\!\!q^2
&A_{4,1,3,2}\!\!=\!\!q\!\!+\!\!q^2\cr
&B_{4,1,4,0}\!\!=\!\!q&A_{4,1,4,1}\!\!=\!\!q\cr
&B_{4,2,1,3}\!\!=\!\!q^5
&A_{4,2,1,4}\!\!=\!\!q^5\cr
&B_{4,2,2,2}\!\!=\!\!q^4\!\!+\!\!q^5
&A_{4,2,2,3}\!\!=\!\!q^4\!\!+\!\!q^5\cr
&B_{4,2,3,1}\!\!=\!\!q^3\!\!+\!\!2q^4\!\!+\!\!q^5
&A_{4,2,3,2}\!\!=\!\!q^3\!\!+\!\!2q^4\!\!+\!\!q^5\cr
&B_{4,2,4,0}\!\!=\!\!2q^3\!\!+\!\!2q^4
&A_{4,2,4,1}\!\!=\!\!2q^3\!\!+\!\!2q^4\cr
&B_{4,3,4,0}\!\!=\!\!q^6
&A_{4,3,4,1}\!\!=\!\!q^6\cr
\noalign{\hrule}
\noalign{\smallskip}
}}}

}
\nobreak
\centerline{Table 2. Polynomials $A_{n,k,a,b}(q)$ and
$B_{n,k,a,b}(q)$ for $0\le n\le 4$}
\nobreak
\centerline{($A_{n,k,a,b}:=A_{n,k,a,b}(q)$ and
$B_{n,k,a,b}:=B_{n,k,a,b}(q)$)}

\vfill\eject
\bigskip\bigskip

{\eightpoint
\centerline{\vtop{\offinterlineskip
\halign{\vrule\hbox{\vrule height 9pt depth 4pt width0pt}
\hfil$#$\hfil\thinspace \vrule &\thinspace \hfil$#$\hfil\thinspace \vrule
&\thinspace \hfil$#$\hfil\thinspace \vrule
&\thinspace \hfil$#$\hfil\thinspace \vrule
&\thinspace \hfil$#$\hfil\thinspace \vrule
\cr
\noalign{\hrule}
m\!=\!&0&1&2&3\cr
\noalign{\hrule}
\kern-3pt
\,n\!=\!0\quad&
&A_{0(00)}\!\!=\!\!1&&\cr
\noalign{\hrule}
1&A_{1(1)}\!\!=\!\!1&
&A_{1(010)}\!\!=\!\!1&\cr
\noalign{\hrule}
2&&A_{2(20)}\!\!=\!\!1
&
&A_{2(010)}\!\!=\!\!1\!\!+\!\!q\cr
&&A_{2(02)}\!\!=\!\!q&&\cr
\noalign{\hrule}
3&A_{3(3)}\!\!=\!\!q\!\!+\!\!q^2
&&A_{3(210)}\!\!=\!\!1\!\!+\!\!q\!\!+\!\!q^2&
\cr
&&
&A_{3(012)}\!\!=\!\!q\!\!+\!\!q^2\!\!+\!\!q^3&\cr
&&&A_{3(030)}\!\!=\!\!q\!\!+\!\!q^2&
\cr
\noalign{\hrule}
4&&A_{4(04)}
\!\!=\!\!q^2\!\!+\!\!q^3\!\!+\!\!2q^4\!\!+\!\!q^5&
&A_{4(0112)}\!\!
=\!\!(1\!\!
+\!\!q\!\!+\!\!q^2)(q\!\!+\!\!q^2\!\!+\!\!q^3
\!\!+\!\!q^4)\cr
&&A_{4(22)}
\!\!=(q\!\!+\!\!q^3)(1\!\!+\!\!q\!\!+\!\!q^2)&
&A_{4(0130)}\!\!=\!\!(q\!\!
+\!\!q^2)(1\!\!+\!\!q\!\!+\!\!q^2+\!\!q^3)\cr
&&A_{4(40)}
\!\!=\!\!q\!\!+\!\!2q^2\!\!+\!\!q^3\!\!+\!\!q^4&
&A_{4(0310)}\!\!=\!\!(q\!\!+\!\!q^2)(1\!\!
+\!\!q\!\!+\!\!q^2
\!\!+\!\!q^3)\cr
&&&&A_{4(2110)}\!\!
=\!\!(1\!\!
+\!\!q\!\!+\!\!q^2)(1\!\!+\!\!q\!\!+\!\!q^2
\!\!+\!\!q^3)\cr
\noalign{\hrule}
\noalign{\medskip}
}}}

\bigskip
\centerline{\vtop{\offinterlineskip
\halign{\vrule\hbox{\vrule height 9pt depth 4pt width0pt}
\hfil$#$\hfil\thinspace \vrule 
&\thinspace \hfil$#$\hfil\thinspace \vrule
&\thinspace \hfil$#$\hfil\thinspace \vrule
\cr
\noalign{\hrule}
m\!=\!&4&5\cr
\noalign{\hrule}
n\!=\!3\quad
&A_{3(01110)}\!\!=\!\!(1\!\!
+\!\!q)(1\!\!+\!\!q\!\!+\!\!q^2)&\cr
\noalign{\hrule}
\ 4&&A_{4(011110)}\!\!=
\!\!(1\!\!
+\!\!q)(1\!\!+\!\!q\!\!+\!\!q^2)
(1\!\!+\!\!q\!\!+\!\!q^2\!\!+\!\!q^3)\cr
\noalign{\hrule}
\noalign{\bigskip}
}}}
}
\centerline{Table 3. Polynomials $A_{n,{\bf c}}(q)$ for $0\le n\le 4$,
($A_{n(c_0\cdots
c_m)}:=A_{n,(c_0,\ldots,c_m)}(q)$)}

{\eightpoint
\bigskip\bigskip\bigskip
\centerline{\vtop{\offinterlineskip
\halign{\vrule\hbox{\vrule height 9pt depth 4pt width0pt}
\hfil$#$\hfil\thinspace \vrule 
&\thinspace \hfil$#$\hfil\thinspace \vrule
&\thinspace \hfil$#$\hfil\thinspace \vrule
&\thinspace \hfil$#$\hfil\thinspace \vrule
&\thinspace \hfil$#$\hfil\thinspace \vrule
&\thinspace \hfil$#$\hfil\thinspace \vrule
&\thinspace \hfil$#$\hfil\thinspace \vrule
\cr
\noalign{\hrule}
m\!=\!&0&1&2&3&4&5\cr
\noalign{\hrule}
n\!=\!0\quad
&1&1&&&&\cr
\noalign{\hrule}
\ 1&1&1&1&&&\cr
\noalign{\hrule}
\ 2&1&1\!\!+\!\!q&1\!\!+\!\!q&1\!\!+\!\!q&&\cr
\noalign{\hrule}
\ 3&q\!\!+\!\!q^2&1\!\!+\!\!2q\!\!+\!\!2q^2
&1\!\!+\!\!3q\!\!+\!\!3q^2\!\!+\!\!q^3
&1\!\!+\!\!2q\!\!+\!\!2q^2\!\!+\!\!q^3
&1\!\!+\!\!2q\!\!+\!\!2q^2\!\!+\!\!q^3&\cr
\noalign{\hrule}
\ 4&q\!\!+\!\!2q^2\!\!+\!\!q^3\!\!+\!\!q^4
&2q\!\!+\!\!4q^2\!\!+\!\!4q^3
&1\!\!+\!\!4q\!\!+\!\!7q^2\!\!+\!\!7q^3
&1\!\!+\!\!5q\!\!+\!\!9q^2\!\!+\!\!1\!0q^3
&1\!\!+\!\!3q\!\!+\!\!5q^2\!\!+\!\!6q^3
&1\!\!+\!\!3q\!\!+\!\!5q^2\!\!+\!\!6q^3\cr
&&\kern5pt{}+\!\!4q^4
\!\!+\!\!2q^5
&\kern13pt{}+\!\!6q^4
\!\!+\!\!3q^5
&\kern5pt{}+\!\!9q^4
\!\!+\!\!5q^5\!\!+\!\!q^6
&\kern5pt{}+\!\!5q^4
\!\!+\!\!3q^5\!\!+\!\!q^6
&\kern5pt{}+\!\!5q^4
\!\!+\!\!3q^5\!\!+\!\!q^6\cr
\noalign{\hrule}
\noalign{\bigskip}
}}}
}
\centerline{Table 4. Polynomials $A_{n,a+b=m}(q)$ 
$(n\equiv m+1\kern-5pt \pmod2)$ and}
\centerline{
$B_{n,a+b=m}(q)$ $(n\equiv m\kern-5pt \pmod2)$ for $0\le
n\le 4$}

\bigskip\bigskip

\bigskip
\vfill\eject

\vglue 25pt
\centerline{\bf References}

{\eightpoint
\bigskip
\article An79|Andr\'e, D\'esir\'e|D\'eveloppement de
$\sec x$ et ${\rm tg}\,x$|C.R. Acad.
Sci. Paris|88|1879|965--967|

\article An81|Andr\'e, D\'esir\'e|Sur les permutations
altern\'ees|J. Math. Pures et Appl.|7|1881|167--184|

\livre An76|Andrews, George E|The Theory of
Partitions|Addison-Wesley, Reading MA, {\oldstyle 1976}
({\it Encyclopedia of Math. and its Appl.}~{\bf 2})|

\livre  AAR00|Andrews, George E.; Askey, Richard; Roy
R|Special Functions|Cambridge University Press, {\oldstyle
2000}|

\article AF80|Andrews, George E.; Foata, Dominique|Congruences
for the $q$-secant number|Europ. J. Combin.|1|1980|283--287|

\article AG78|Andrews, George E.; Gessel,
Ira|Divisibility properties
of the $q$-tangent numbers|Proc.
Amer. Math. Soc.|68|1978|380--384|

\article Ar92|Arnold, V. I|Springer numbers and
Morsification spaces|J. Algebraic Geom.|1|1992|197--214|

\article Ar92a|Arnold, V. I|The calculus of snakes 
and the combinatorics of Bernoulli, Euler and Springer
numbers of Coxeter groups|Uspekhi Mat. nauk.|47|1992|3--45
= Russian Math. Surveys, {\bf47} ({\oldstyle1992}), 1--51|

\article Ca54|Carlitz, Leonard|$q$-Bernoulli and 
Eulerian numbers|Trans.
Amer. Math. Soc.|76|1954|332--350|

\article Ca75|Carlitz, Leonard|A combinatorial property of
$q$-Eulerian numbers|Amer. Math. Monthly|82|1975|51--54|

\article CS72|Carlitz, Leonard; Scoville, Richard|Tangent 
Numbers and Inversions|Duke Math. J.|39|1972|413--429|

\livre Co74|Comtet, Louis|Advanced Combinatorics|D.
Reidel/Dordrecht-Holland, Boston, {\oldstyle 1974}|

\divers Eu48|Euler, Leonhard|De partitione numerorum 
{\it in} Introductio in analysin infinitorum,
	chap. 16. Opera Omnia, I8, {\oldstyle 1748},
	p. 313--338. B.G. Teubner [Adolf Krazer, Ferdinand Rudio,
	eds.], {\oldstyle 1922}|

\article Du95|Dumont, Dominique|Further triangles of 
Seidel-Arnold type and continued fractions related to 
Euler and Springer numbers|Adv. Appl. Mah.|16|1995|275--296|

\divers Du12|Dumont, Dominique|Private communication, {\oldstyle 2012}|

\article Fo68|Foata, Dominique|On the Netto inversion 
number of a sequence|Proc. Amer. Math. Soc.|19|1968|236--240|

\article Fo81|Foata, Dominique|Further divisibility properties of
the $q$-tangent numbers|Proc. Amer. Math. Soc.|81|1981|143--148|

\divers FH11|Foata, Dominique; Han, Guo-Niu|The
$(t,q)$-Analogs of Secant and Tangent Numbers, {\sl
Electron. J. Combin.} {\bf 18} ({\oldstyle 2011}) (The
Zeilberger Festschrift), \#P7, 16 p|

\article FS78|Foata, Dominique; Sch\"utzenberger,
Marcel-Paul|Major Index and Inversion number of
Permutations|Math. Nachr.|83|1978|143--159|

\article Fr78|Fran\c con, Jean|Histoires de fichiers|RAIRO 
Inf. Theor.|12|1978|49--62|

\livre GR90|Gasper, George; Rahman, Mizan|Basic hypergeometric
series| Encyclopedia of Math. and its Appl. {\bf 35},
Cambridge Univ. Press, Cambridge, {\oldstyle 1990}|

\article Gl98|Glaisher, J. W. L|On the Bernoullian
function|Quart. J. Pure Appl. Math.|29|1898|1--168|

\article Gl99|Glaisher, J. W. L|On a set of
coefficients analogous  to the Eulerian numbers|Proc.
London Math. Soc.|31|1899|216--235|

\article Gl14|Glaisher, J. W. L|On the coefficients 
in the expansions of $\cos x/\!\cos 2x$ and $\sin x/\!\cos
2x$|Quart. J. Pure Appl. Math.|45|1914|187--222|

\article Ho95|Hoffman, Michael E|Derivative Polynomials for
Tangent and Secant|Amer. Math. Monthly|102|1995|23--30|

\divers Ho99|Hoffman, Michael E|Derivative polynomials, Euler
polynomials, and associated integer sequences, {\sl Electron.
J. Combin.} {\bf 6} ({\oldstyle 1999}), \#R21|

\article Ja04|Jackson, F.H|A basic-sine and cosine with
symbolic solutions of certain differential equations|Proc.
Edinburgh Math. Soc.|22|1904|28--39|

\divers Jo11|Josuat-Verg\`es, Matthieu|Enumeration of
snakes and cycle-alternating
permutations, arXiv:1011.0929,  {\oldstyle 2011}| 

\divers JNT12|Josuat-Verg\`es, Matthieu; Novelli, Jean-Chrisophe;
Thibon, Jean-Yves|The algebraic combinatorics of snakes,
arXiv:1110.5272, {\oldstyle 2012}|

\article KB67|Knuth, D. E.; Buckholtz, Thomas J|Computation
of tangent, Euler and Bernoulli numbers|Math.
Comp.|21|1967|663--688|

\livre Lo83|Lothaire, M|Combinatorics on
Words|Addison-Wesley Publ. Co., Reading, Mass., {\oldstyle 1983} (Encyclopedia
		of Math. and its Appl., {\bf 17})|

\livre Ni23|Nielsen, Niels|Trait\'e \'el\'ementaire des nombres
de Bernoulli|Paris, Gau\-thier-Villars, {\oldstyle 1923}|

\article Sp71|Springer, T.A|Remarks on a combinatorial
problem|Nieuw Arch. Wisk.|19|1971|30--36|

\divers Sl06|Sloane, Neil J. A|The On-Line Encyclopedia of
Integer Sequences (OEIS), {\oldstyle 2006},
{\tt http://oeis.org/}|

\article St76|Stanley, Richard P|Binomial posets, M\"obius inversion, and permutation enumeration|J. Combin. Theory Ser. A|20|1976|336--356|

\livre St97|Stanley, Richard P|Enumerative
Combinatorics|Vol. 1, Cambridge University Press,
{\oldstyle 1997}|

\article St10|Stanley, Richard P|A Survey of Alternating Permutations|Contemporary Mathematics|531|2010|165-196|

}
\bigskip
\hbox{\qquad\qquad\vtop{\halign{#\hfil\cr
Dominique Foata \cr
Institut Lothaire\cr
1, rue Murner\cr
F-67000 Strasbourg, France\cr
\noalign{\smallskip}
{\tt foata@unistra.fr}\cr}}
\qquad
\vtop{\halign{#\hfil\cr
Guo-Niu Han\cr
I.R.M.A. UMR 7501\cr
Universit\'e de Strasbourg et CNRS\cr
7, rue Ren\'e-Descartes\cr
F-67084 Strasbourg, France\cr
\noalign{\smallskip}
{\tt guoniu.han@unistra.fr}\cr}}
}

\bye